\documentclass[a4paper]{amsart}
\usepackage{etex}
\usepackage[OT2, T1]{fontenc}
\usepackage{url}
\usepackage{amsmath}
\usepackage{lmodern}
\usepackage{array}
\usepackage{graphicx}
\usepackage{amsfonts}
\usepackage{amssymb}

\usepackage{xargs}
\usepackage[pdftex, dvipsnames]{xcolor}
\usepackage[colorinlistoftodos,prependcaption,textsize=tiny]{todonotes}
\newcommandx{\unsure}[2][1=]{\todo[linecolor=red,backgroundcolor=red!25,bordercolor=red,#1]{#2}}
\newcommandx{\change}[2][1=]{\todo[linecolor=blue,backgroundcolor=blue!25,bordercolor=blue,#1]{#2}}
\newcommandx{\info}[2][1=]{\todo[linecolor=OliveGreen,backgroundcolor=OliveGreen!25,bordercolor=OliveGreen,#1]{#2}}
\newcommandx{\improvement}[2][1=]{\todo[linecolor=Plum,backgroundcolor=Plum!25,bordercolor=Plum,#1]{#2}}
\newcommandx{\thiswillnotshow}[2][1=]{\todo[disable,#1]{#2}}

\definecolor{imperialred}{RGB}{237, 41, 57}
\definecolor{royalblue}{RGB}{64, 106, 212}
\definecolor{link}{RGB}{11,0,128}
\usepackage[colorlinks=true, citecolor=royalblue,linkcolor=imperialred,urlcolor=link]{hyperref}
\usepackage{amsthm}
\usepackage{esint}

\usepackage{paralist}
\usepackage{colonequals}		
\usepackage{mathabx}		
\usepackage[left]{lineno}
\usepackage{amscd}
\usepackage[all,cmtip]{xy}	
\usepackage{sseq}			
\usepackage{verbatim}
\usepackage{parskip}
\usepackage{microtype}
\usepackage[in]{fullpage}
\usepackage{graphicx}
\usepackage{mathrsfs} 			
\usepackage{bm} 				
\usepackage{cleveref}
\usepackage{enumitem}
\usepackage{moreenum}			
\usepackage{stmaryrd}
\usepackage{blindtext}
\usepackage{tikz-cd}
\usepackage{tikz}
\usetikzlibrary{bending,matrix,arrows,decorations.pathmorphing,backgrounds,positioning,fit,petri,arrows.meta}
\usepackage{textcomp}
\usepackage{indentfirst}
\usetikzlibrary{arrows}
\usepackage{enumitem}
\tikzset{commutative diagrams/.cd,arrow style=tikz,diagrams={>=latex'}}

\usepackage{filecontents}
\usepackage[alphabetic,lite]{amsrefs} 	
\usepackage{stmaryrd}	
\usepackage[foot]{amsaddr}
\usepackage{subfiles}
\usepackage{ragged2e}

\DeclareSymbolFont{cyrletters}{OT2}{wncyr}{m}{n}
\DeclareMathSymbol{\Sha}{\mathalpha}{cyrletters}{"58}

\newcommand{\GG}{\Gamma}

\newcommand{\GF}{\Phi}

\newcommand{\bA}{\mathbb{A}}

\newcommand{\bG}{\mathbb{G}}

\newcommand{\cT}{\mathcal{T}}

\newcommand{\cV}{\mathcal{V}}

\newcommand{\cX}{\mathcal{X}}

\newcommand{\fm}{\mathfrak{m}}

\newcommand{\sG}{\mathscr{G}}

\newcommand{\sI}{\mathscr{I}}

\newcommand{\ra}{\rightarrow}

\newcommand{\wt}{\widetilde}
\newcommand{\wh}{\widehat}

\newcommand{\pr}{^{\prime}}
\newcommand{\prpr}{^{\prime\prime}}
\newcommand{\op}{^{\mathrm{op}}}		
\newcommand{\ce}{\colonequals}
\newcommand{\ov}{\overline}
\newcommand{\un}{\underline}

\renewcommand{\b}{\textbf}
\newcommand{\surjects}{\twoheadrightarrow}
\newcommand{\injects}{\hookrightarrow}
\newcommand{\isom}{\simeq}			

\newcommand{\isomto}{\overset{\sim}{\longrightarrow}}

\newcommand{\et}{\mathrm{\acute{e}t}}	

\makeatletter
\DeclareRobustCommand\widecheck[1]{{\mathpalette\@widecheck{#1}}}
\def\@widecheck#1#2{%
\setbox\z@\hbox{\m@th$#1#2$}%
\setbox\tw@\hbox{\m@th$#1%
\widehat{%
\vrule\@width\z@\@height\ht\z@
\vrule\@height\z@\@width\wd\z@}$}%
\dp\tw@-\ht\z@
\@tempdima\ht\z@ \advance\@tempdima2\ht\tw@ \divide\@tempdima\thr@@
\setbox\tw@\hbox{%
\raise\@tempdima\hbox{\scalebox{1}[-1]{\lower\@tempdima\box
\tw@}}}%
{\ooalign{\box\tw@ \cr \box\z@}}}
\makeatother

\newcommand{\leftsub}[2]{{\vphantom{#2}}_{#1}{#2}}
\providecommand{\SP}[1]{\cite[\href{https://stacks.math.columbia.edu/tag/#1}{#1}]{SP}}

\providecommand{\p}[1]{\left(#1\right)}
\providecommand{\hkh}[1]{\left\{#1\right\}}

\providecommand{\ucolon}{{\upshape:} }
\providecommand{\uscolon}{{\upshape;} }

\providecommand{\fps}[1]{\llbracket#1\rrbracket}
\providecommand{\lps}[1]{(\!(#1)\!)}
\DeclareMathOperator{\Spec}{Spec}		
\DeclareMathOperator{\rad}{rad}			
\DeclareMathOperator{\Frac}{Frac}		
\DeclareMathOperator{\id}{id}			
\DeclareMathOperator{\Gal}{Gal}	
\DeclareMathOperator{\Res}{Res}		
\DeclareMathOperator{\Lie}{Lie}		
\newcommand{\ba}{\begin{aligned}}
\newcommand{\ea}{\end{aligned}}
\newcommand{\be}{\begin{equation}}
\newcommand{\ee}{\end{equation}}
\newcommand{\pf}{\begin{proof}}
\newcommand{\bpf}{\begin{proof}}
\newcommand{\epf}{\end{proof}}
\newcommand{\bthm}{\begin{thm}}
\newcommand{\ethm}{\end{thm}}

\newcommand{\bprop}{\begin{prop}}
\newcommand{\eprop}{\end{prop}}
\newcommand{\bcor}{\begin{cor}}
\newcommand{\ecor}{\end{cor}}

\newcommand{\brem}{\begin{rem}}
\newcommand{\erem}{\end{rem}}

\newcommand{\mtg}{\underline{\mathrm{Tor}}}
\newcommand{\brems}{\begin{rems} \hfill \begin{enumerate}[label=\b{\thenumberingbase.},ref=\thenumberingbase]}
\newcommand{\remi}{\addtocounter{numberingbase}{1} \item}
\newcommand{\erems}{\end{enumerate} \end{rems}}

\newcommand{\begs}{\begin{egs} \hfill \begin{enumerate}[label=\b{\thenumberingbase.},ref=\thenumberingbase]}

\newcommand{\eegs}{\end{enumerate} \end{egs}}

\newcommand{\eremstweak}{\end{enumerate} \end{rems-tweak}}
\newcommand{\eremst}{\end{enumerate} \end{rems-tweak}}
\newcommand{\blem}{\begin{lemma}}
\newcommand{\elem}{\end{lemma}}

\newcommand{\bconj}{\begin{conj}}
\newcommand{\econj}{\end{conj}}
\newcommand{\bprob}{\begin{Problem}}
\newcommand{\eprob}{\end{Problem}}

\newcommand{\bq}{\begin{Q}}
\newcommand{\eq}{\end{Q}}
\newcommand{\benum}{\begin{enumerate}[label={{\upshape(\alph*)}}]}
\newcommand{\benuma}{\begin{enumerate}[label={{\upshape(\arabic*)}}]}
\newcommand{\benumr}{\begin{enumerate}[label={{\upshape(\roman*)}}]}
\newcommand{\eenum}{\end{enumerate}}
\newcommand{\bc}{}
\newcommand{\bd}{\begin{defn}}
\newcommand{\ed}{\end{defn}}

\newcommand{\beg}{\begin{eg}}
\newcommand{\eeg}{\end{eg}}

\newcommand{\bcl}{\begin{claim}}
\newcommand{\ecl}{\end{claim}}

\newcommand{\q}{\quad}
\newcommand{\qq}{\quad\quad}
\newcommand{\qqq}{\quad\quad\quad}
\newcommand{\qqqq}{\quad\quad\quad\quad}

\newcommand{\qqqqqqq}{\quad\quad\quad\quad\quad\quad\quad}
\newcommand{\qqqqqqqq}{\quad\quad\quad\quad\quad\quad\quad\quad}

\newcommand{\tst}{\textstyle}

\newcommand{\outauto}{\un{\mathrm{Out}}}

\usepackage{aliascnt}
\newaliascnt{numberingbase}{subsection}

\theoremstyle{plain}
\newtheorem{thm}[numberingbase]{Theorem}
\Crefname{thm}{Theorem}{Theorems}

\Crefname{rethm}{Theorem}{Theorem}
\newtheorem{prop}[numberingbase]{Proposition}
\Crefname{prop}{Proposition}{Propositions}
\newtheorem{Q}[numberingbase]{Question}
\Crefname{Q}{Question}{Questions}
\newtheorem{Problem}[subsection]{Problem}
\Crefname{Problem}{Problem}{Problems}
\newtheorem{conj}[numberingbase]{Conjecture}
\Crefname{conj}{Conjecture}{Conjectures}
\newtheorem{cor}[numberingbase]{Corollary}
\Crefname{cor}{Corollary}{Corollaries}
\newtheorem{lemma}[numberingbase]{Lemma}

\Crefname{subprop}{Proposition}{Propositions}

\Crefname{subcor}{Corollary}{Corollaries}

\Crefname{sublem}{Lemma}{Lemmas}

\theoremstyle{remark}
\newtheorem{claim}[equation]{Claim}
\Crefname{claim}{Claim}{Claims}

\Crefname{subrem}{Remark}{Remarks}

\theoremstyle{definition}
\newtheorem{defn}[numberingbase]{Definition}
\Crefname{defn}{Definition}{Definitions}

\Crefname{conv}{Convention}{Conventions}
\newtheorem{eg}[numberingbase]{Example}
\Crefname{eg}{Example}{Examples}
\newtheorem{rem}[numberingbase]{Remark}
\Crefname{rem}{Remark}{Remarks}
\newtheorem*{rems}{Remarks}
\newtheorem*{egs}{Examples}

\theoremstyle{plain}
\newtheorem{thm-tweak}[subsection]{Theorem}
\Crefname{thm-tweak}{Theorem}{Theorems}
\newtheorem{lemma-tweak}[subsection]{Lemma}
\Crefname{lemma-tweak}{Lemma}{Lemmas}
\newtheorem{cor-tweak}[subsection]{Corollary}
\Crefname{cor-tweak}{Corollary}{Corollaries}
\newtheorem{prop-tweak}[subsection]{Proposition}
\Crefname{prop-tweak}{Proposition}{Propositions}
\newtheorem{conj-tweak}[subsection]{Conjecture}
\Crefname{conj-tweak}{Conjecture}{Conjectures}

\theoremstyle{definition}
\newtheorem{defn-tweak}[subsection]{Definition}
\Crefname{defn-tweak}{Definition}{Definitions}
\newtheorem{eg-tweak}[subsection]{Example}
\Crefname{eg-tweak}{Example}{Examples}
\newtheorem*{rems-tweak}{Remarks}
\newtheorem{rem-tweak}[subsection]{Remark}
\Crefname{rem-tweak}{Remark}{Remarks}

\newtheoremstyle{subsection-tweak}
   {11pt}
   {3pt}%
   {}
   {}%
   {\bfseries}
   {}%
   {.5em}
   {\thmnumber{\@{#1}{}\@{#2}.}%
	\thmnote{~{\bfseries#3.}}}

\theoremstyle{subsection-tweak}
\newtheorem{pp}[numberingbase]{}
\newcommand{\bpp}{\begin{pp}}
\newcommand{\epp}{\end{pp}}

\theoremstyle{subsection-tweak}
\newtheorem{pp-tweak}[subsection]{}

\numberwithin{equation}{numberingbase}

\makeatletter
\def\@tocline#1#2#3#4#5#6#7{
	\begingroup 
	\@ifempty{#4}{%
	}{%
	}%

	\parindent\z@ \leftskip#3\relax \advance\leftskip\@tempdima\relax
	#5\hskip-\@tempdima
	  \ifcase #1
	   \or\or \hskip 2em \or \hskip 1em \else \hskip 3em \fi%
	  #6\nobreak\relax
	\dotfill\hbox to\@pnumwidth{\@tocpagenum{#7}}\par
	\nobreak
	\endgroup
  }
 \def\l@section{\@tocline{1}{0pt}{1pc}{}{}}

\renewcommand{\tocsection}[3]{%
  \indentlabel{\@ifnotempty{#2}{\makebox[1.3em][l]{%
	\ignorespaces#1 \bfseries{#2}.\hfill}}}\bfseries{#3}
	\vspace{-0.5pt}}

\renewcommand{\tocsubsection}[3]{%
  \indentlabel{\@ifnotempty{#2}{\hspace*{-0.5em}\makebox[2.1em][l]{%
	\ignorespaces#1#2.\hfill}}}#3
	\vspace{1.5pt}}

\makeatother

\makeatletter
\newcommand\appendix@section[1]{%
  \refstepcounter{section}%
  \orig@section*{Appendix \@Alph\c@section. #1}%
 \addcontentsline{toc}{section}{Appendix \@Alph\c@section. #1}%
}
\let\orig@section\section
\g@addto@macro\appendix{\let\section\appendix@section}
\makeatother

\DeclareMathVersion{normal2}

\begin{document}

\title{THE GROTHENDIECK--SERRE CONJECTURE OVER SEMILOCAL DEDEKIND RINGS}

\author{NING GUO}
\address{Laboratoire de Math\'ematiques d'Orsay, Univ. Paris-Sud, Universit\'e Paris-Saclay, 91405 Orsay, France}
\email{ning.guo@math.u-psud.fr}

\begin{abstract}
For a reductive group scheme $G$ over a semilocal Dedekind ring $R$ with total ring of fractions $K$, we prove that no nontrivial $G$-torsor trivializes over $K$. This generalizes a result of Nisnevich--Tits, who settled the case when $R$ is local. Their result, in turn, is a special case of a conjecture of Grothendieck--Serre that predicts the same over any regular local ring. With a patching technique and weak approximation in the style of Harder, we reduce to the case when $R$ is a complete discrete valuation ring. Afterwards, we consider Levi subgroups to reduce to the case when $G$ is semisimple and anisotropic, in which case we take advantage of Bruhat--Tits theory to conclude. Finally, we show that the Grothendieck--Serre conjecture implies that any reductive group over the total ring of fractions of a regular semilocal ring $S$ has at most one reductive $S$-model.
\end{abstract}
 \maketitle

\hypersetup{
    linktoc=page,     
}
\renewcommand*\contentsname{}
\q\\
\tableofcontents
\section{Introduction}\label{intro}
The Grothendieck--Serre conjecture was proposed by J.-P. Serre (\cite{Ser58}*{p.~31, Rem.}) and A. Grothendieck (\cite{Gro58}*{pp.~26-27, Rem.~3}) in 1958, who predicted that for an algebraic group $G$ over an algebraically closed field, a $G$-torsor over a nonsingular variety is Zariski-locally trivial if it is generically trivial. Subsequently, Grothendieck extended the conjecture to semisimple group schemes over regular schemes (\cite{Gro68}*{Rem.~1.11.a}). By spreading out, the conjecture reduces to its local case whose precise statement is the following.
		\bconj[Grothendieck--Serre] \label{GSconjgeneral}
For a reductive group $G$ over a regular local ring $R$ with fraction field $K$, a $G$-torsor that becomes trivial over $K$ is trivial. In other words, the following map between nonabelian \'etale cohomology pointed sets has trivial kernel\ucolon 
\be\label{gsmap}
H^1_{\et}(R,G)\ra H^1_{\et}(K,G).
\ee
\econj
In this paper, we prove a variant of \Cref{GSconjgeneral} when $R$ is a semilocal Dedekind ring.
Our main result is the following theorem.
\bthm \label{GSconj}
  For a reductive group scheme $G$ over a semilocal Dedekind ring $R$ with total ring of fractions $K$, the following pullback map of nonabelian \'etale cohomology pointed sets
  \[
   \text{$H^1_{\et}(R, G)\hookrightarrow H^1_{\et}(K, G)$ \qq is injective.}
  \]
\ethm

Beyond the trivial case of $0$-dimensional $R$, \Cref{GSconjgeneral} has several known cases.
\begin{itemize}
  \item[(i)] The case when $R$ is a discrete valuation ring with a perfect residue field was addressed in Nisnevich's PhD thesis \cite{nisnevich1983etale}*{2, Thm.~7.1} and his Comptes Rendus paper \cite{Nis84}*{Thm.~4.2}. There, he reduced to the complete case that had been considered in unpublished work of Tits that rested on Bruhat--Tits theory (see \cite{BrT3}*{Lem.~3.9}, but still with auxiliary conditions). 
  Various other special cases are based on the discrete valuation ring case, for instance, the case when $R$ is of arbitrary dimension and complete regular local, see \cite{CTS79}*{6.6.1}; further, the case when $R$ is Henselian and $G$ splits was proved explicitly in \cite{BB70}*{Prop.~2}.
  \item[(ii)] The case when $\dim R=2$ with infinite residue field and $G$ is quasi-split was considered by Nisnevich in \cite{Nis89}*{Thm.~6.3}.
  \item[(iii)] The case when $G$ is torus was settled by Colliot-Th\'el\`ene and Sansuc in \cite{CTS87}*{Thm.~4.1}. This result is often useful for various reductions of more general cases.
  \item[(iv)] The case when $R$ is a semilocal Dedekind domain and $G$ is a form of $PGL_n$, $PSp_n$ and $PO_n$ ($PO_n$ is not connected) as a variant of the Grothendieck--Serre conjecture was proved by Beke and Van Geel in \cite{BVG14}*{Thm.~3.7}.
  \item[(v)] The case when $R$ contains a field $k$ is summarized as follows. 
  When $k$ is algebraically closed and $G$ is defined over $k$, the conjecture was settled by Colliot-Th\'el\`ene--Ojanguren in \cite{CTO92}.
  For regular semilocal domains $R$ containing a field $k$, when $k$ is infinite, the conjecture was recently proved by Fedorov and Panin in \cite{FP15}; the case when $k$ is finite was originally proved by Panin in \cite{Pan15}, whose exposition was organized in his latter work \cite{Pan17}. 
  For a certain simplification of \emph{ibid.}, see \cite{FR18}. 
  \item[(vi)] The case when $R$ is a semilocal Dedekind domain and $G$ is a form of $GL_n$, $O_n$ and $Sp_n$ was recently settled by Bayer-Fluckiger and First in \cite{BF17}*{Thm.~5.3}, where they also proved a variant for non-reductive groups with generic fibers of the form $GL_n$, $O_n$ and $Sp_n$. They also investigated the corresponding problem for $PGL_n$, $PO_n$ and $PSp_n$ in \cite{BF18}, where counterexamples for non-reductive $G$ emerge, even when $R$ is a complete discrete valuation ring.
  \item[(vii)] The case when $R$ is a semilocal Dedekind domain and $G$ is a semisimple simply connected $R$-group scheme such that every semisimple normal $R$-subgroup scheme of $G$ is isotropic was recently proved by Panin and Stavrova in \cite{PS17}*{Thm.~3.4}. By induction on the number of maximal ideals of $R$ and a decomposition of groups, they reduce to the case when $R$ is a Henselian discrete valuation ring and use Nisnevich's result \cite{Nis84}*{Thm.~4.2} to conclude.
  \item[(viii)] A stronger conjecture predicts that the map (\ref{gsmap}) is injective, but it is equivalent to \Cref{GSconjgeneral} due to a twisting technique (\Cref{GSconjtrivialkernel}), which is standard.
\end{itemize}

The main result \Cref{GSconj} covers (i) and the corresponding cases of (iv)---(vii) and finishes the remaining cases of the semilocal Dedekind variant of \Cref{GSconjgeneral}.
To be precise, this article is aimed at working out Nisnevich's proof \cite{nisnevich1983etale}, extending it to cover $1$-dimensional regular semilocal rings, and in particular, generalizing the result of Panin--Stavrova \cite{PS17}*{Thm.~3.4} by eliminating the semisimple simply connected isotropic condition on the reductive group scheme $G$.

We now summarize Nisnevich's work (\cite{nisnevich1983etale},\cite{Nis84}) in the discrete valuation ring case. 
For a reductive group scheme $G$ over a discrete valuation ring $R$ with fraction field $K$, we let $\wh{R}$ be the completion of $R$ and $\wh{K}\ce \Frac(\wh{R})$.
 Built upon Harder's weak approximation \cite{Harder}, Nisnevich proved that $G(\wh{K})$ decomposes as a product of $G(\wh{R})$ and $G(K)$. This reduces one to the case when $R$ is complete. Following an argument of Tits in this case, one uses Bruhat--Tits theory and Galois cohomology to conclude.

Unfortunately, the details of Tits' argument, which are given in \cite{nisnevich1983etale}, have some unclear points.
These details relied on \cite{BrT2}, which appeared in print two years after \cite{nisnevich1983etale}.
Further, Nisnevich's proof is written under an auxiliary hypothesis that the residue field of $R$ is perfect (\cite{nisnevich1983etale}*{Ch.~2, Thm.~7.1}). Although this hypothesis does not appear in \cite{Nis84}, the part of the proof corresponding to the key step \cite{nisnevich1983etale}*{Ch.~2, 4.2} is implicit as an unpublished result of Tits.  Moreover, 
in \cite{nisnevich1983etale}*{Ch.~2, Thm.~4.2}, there is a decomposition of groups, which fails in general and forces us to consider alternative reductions (see Remark \ref{Nisdefectproduct}). Besides, some recent articles such as \cite{PS17}*{Thm.~3.4} use Nisnevich's result \cite{Nis84}*{Thm.~4.2} directly without treating these gaps.

In this article, we first review weak approximation in \Cref{weak approximation}.
For a semilocal Dedekind domain $R$ with fraction field $K$, we consider a reductive $K$-group $G$ instead of merely semisimple as in Harder's original setting in \cite{Harder}.   The goal of \Cref{weak approximation} is to construct an open normal subgroup of the closure of $G(K)$ in $\prod_{v} G(K_v)$ (\Cref{opennormalsubgroup}), where $v$ ranges over the nongeneric points of $\Spec R$ and $K_v$ is the completion of $K$ at $v$. In \Cref{decomposition of groups}, with the aid of Harder's construction, we exhibit the decomposition
\[
\tst \text{\qqqqqqq\qqqqqqqq $\prod_{v} G(K_v)=G(K)\prod_{v} G(R_v)$\qqq \qqqqqqq\q (\Cref{grpdecomp})}
\]
by mildly simplifying the argument in \cite{nisnevich1983etale}*{Ch.~2, \S6}. This permits us to reduce \Cref{GSconj} to the case of complete discrete valuation rings by a  patching technique (\Cref{reducetocompletecase}).

To improve and extend Nisnevich's proof, after reducing to the complete case, we take a different approach. The passage from $G$ to its minimal parabolic subgroup $P$ and the isomorphism $H^1(R, P)\isom H^1(R, L)$ for a Levi subgroup $L$ of $P$ facilitate reduction to the case when $G$ is semisimple and anisotropic (\Cref{reducetosemisimpleanisotropic}). This reduction and subsequent steps rely on properties of anisotropic groups described in \Cref{anisotropicrationalintegerpoints},  where we use the formalism of Bruhat--Tits theory from \cite{BrT1, BrT2}.
\bprop[\Cref{anisotropicrationalintegerpoints}]
  For a reductive group $G$ over a Henselian discrete valuation ring $R$,
  \benum
  \item for the strict Henselization $\wt{R}$ of $R$ with fraction field $\wt{K}$, the group $G(\wt{R})$ is a maximal parahoric subgroup of $G(\wt{K})$;
  \item if $G$ is $K$-anisotropic, then $G(K)=G(R)$.
  \eenum
\eprop
In the final step, we establish the case when $G$ is semisimple anisotropic and $R$ is a complete discrete valuation ring (\Cref{finalproof}) by following the argument of Nisnevich--Tits.

Finally, for a reductive group $G$ over the function field $K$ of a scheme $X$, we consider the number of reductive $X$-models of $G$. In fact, if a variant of \Cref{GSconjgeneral} holds for regular semilocal $X$, then we prove that such models should be unique:
\bprop[\Cref{Uniqueness}]\label{uniq}
For a regular semilocal ring $S$ with total ring of fractions $K$, if
\[
(\ast)\qq \text{for each reductive group scheme $G'$ over $S$, the map $H^1_{\et}(S, G')\ra H^1_{\et}(K, G')$ is injective},\]
then any reductive $K$-group $G$ has at most one reductive $S$-model. 
\eprop
In particular, when $S$ contains a field, the condition $(\ast)$ is satisfied thanks to the corresponding settled case of the Grothendieck--Serre conjecture \cite{FP15, Pan17}, hence follows the uniqueness of reductive models over regular semilocal rings containing a field.
On the other hand, the combination of \Cref{uniq} and our main result \Cref{GSconj} implies the uniqueness of reductive models for semilocal Dedekind rings (\Cref{modelsemilocalDedekind}), which generalizes \cite{Nis84}*{Thm.~5.1} where $S$ is a Henselian local ring and $G$ is semisimple. 

\bpp[Notation and conventions]\label{convention}
If not particularly indicated, in the sequel, we let $R$ be a semilocal Dedekind ring and let $K$ be its total ring of fractions. The completion of $R$ at its radical ideal is a direct product of fields and complete discrete valuation rings $R_v$ with fraction fields $K_v$. In the case when $R$ is local, we also denote the strict Henselization of $R$ by $\wt{R}$ and its fraction field by $\wt{K}$.
\epp

We also assume that $G$ is a reductive group scheme over $R$ (or over $K$ in \Cref{weak approximation}), that is, $G$ is a smooth, affine $R$-group scheme with connected reductive algebraic groups as geometric fibers. For an $R$-algebra $R'$, the reductive group $G$ is called $R'$\emph{-anisotropic}, if $G_{R'}$ contains no copy of $\mathbb{G}_{m,R'}$.

\subsection*{Acknowledgements} I thank my advisor K\k{e}stutis \v{C}esnavi\v{c}ius. His guidance and advice were helpful. He also provided exhaustive comments for revising. I learnt much from him in several aspects. I also thank Uriya First for giving me useful comments and advice. I thank referees for their careful reading and informative suggestions.

\section{Construction of open normal subgroups} \label{weak approximation}
Let $H$ be a reductive group over a field $F$ equipped with a finite set $\cV$ of nonequivalent nontrivial valuations of rank one.
This section is devoted to working out Harder's construction from \cite{Harder} of an open normal subgroup $N$ contained in the closure of $H(F)$ in $\prod_{v\in \cV}H(F_v)$. Here, $F_v$ is the completion of $F$ at $v$; the groups $H(F_v)$ are endowed with their $v$-adic topologies, and $\prod_{v\in \cV}H(F_v)$ has the product topology. 
In the sequel, we regard $H(F_v)$ as a subgroup of $\prod_{v\in \cV}H(F_v)$ by identifying it with the subgroup $H(F_v)\times \prod_{v\pr\neq v}\{1\}$.
We will need the group $N$ in \Cref{decomposition of groups} when exhibiting the decomposition in \Cref{grpdecomp}, which leads to the reduction of \Cref{GSconj} to the case of complete discrete valuation rings (\Cref{reducetocompletecase}). 
To begin, we recall Grothendieck's Theorem (\cite{SGA3II}*{XIV, 1.1}) that any smooth group of finite type over a field contains a maximal torus. Let $L_w$ be a minimal splitting field of a maximal $F_v$-torus $T$ of $H_{F_v}$, where $w$ is a valuation extending $v$. 
It turns out that the image $U$ of the norm map
\[
N_{L_w/F_v}: T(L_w)\ra T(F_v)
\]
is an open subgroup of $T(F_v)$ and contained in $\ov{H(F)}\cap H(F_v)$ (see \Cref{maxtorusopensubset}), where the closure $\ov{H(F)}$ is formed in $\prod_{v\in \cV}H(F_v)$.
 We use $U$ to construct an open normal subgroup of $H(F_v)$ contained in $\ov{H(F)}\cap H(F_v)$. This gives rise to the desired open normal subgroup $N$ constructed in \Cref{opennormalsubgroup}.
\blem \label{maxtorusopensubset}
For a maximal torus $T$ of $H_{F_v}$, the image $U$ of the norm map
\[
N_{L_w/F_v}\colon T(L_w)\ra T(F_v)
\]
is an open subgroup of $T(F_v)$ contained in $\ov{H(F)}\cap H(F_v)$.
\elem
\bpf
We recall from \cite{SGA3II}*{XIV, 6.1} that the functor
\[
\mtg(H)\colon \mathbf{Sch}\op_{/F}\ra \mathbf{Sets}, \qqq F'\mapsto \{
\text{maximal tori of $H_{F'}$}\}
\]
is representable by an $F$-scheme that is a rational variety (that is, by \cite[\href{https://stacks.math.columbia.edu/tag/0BXP}{0BXP}]{SP}, it has an open dense subscheme isomorphic to an open subscheme of $\mathbb{A}^n_F$). In particular, by the same argument as in \cite{PV}*{Prop.~7.3}, $\mtg(H)$ satisfies weak approximation: $\mtg(H)(F)$ is dense in $\prod_{v\in \cV}\mtg(H)(F_v)$ under the diagonal embedding.

By \cite{SGA3IIInew}*{XXII, 5.8.3}, the morphism
\[
  H_{F_v}\ra \mtg(H_{F_v}), \qq g\mapsto gTg^{-1}
\]
induces an isomorphism $H_{F_v}/\un{\mathrm{Norm}}_{H_{F_v}}(T)\isom \mtg(H_{F_v})$. By \cite{SGA3II}*{XI, 2.4bis}, $\un{\mathrm{Norm}}_{H_{F_v}}(T)$ is smooth. Therefore, by \cite{Ces15d}*{Prop.~4.3}, the map
\[
\phi: H(F_v)\ra \mtg(H)(F_v), \qq g\mapsto gTg^{-1}
\]
is open. Thus, for any neighbourhood $V$ of the neutral element of $H(F_v)$, by weak approximation for $\mtg(H)$, the image $\phi(V)$ contains a maximal torus $T\pr$ of $H$.
Let $L/F$ be a minimal splitting field of $T'$.
Then, we have the decomposition:
\[
\tst L\otimes_F F_v\isom \prod_{i=1}^{r}L_i.
\]
The isomorphism $T'_{F_v}\isom T$ implies that $L_i\isom L_w$ for each $i$. The norm map
$N_{L_w/F_v}: T'(L_w)\ra T'(F_v)$ is the map on the $F_v$-points induced by the Weil restriction $N: \Res_{L_w/F_v}T'_{L_w}\ra T'_{F_v}$, whose kernel is a torus.
So \cite{Ces15d}*{Prop.~4.3} implies the openness of the image $U\pr\ce N_{L_w/F_v}(T\pr(L_w))\subset T\pr(F_v)$.

Now we prove that $U'\subset \ov{H(F)}\cap H(F_v)$. For $a\in T'(L_w)$ and $b\ce N_{L_w/F_v}(a)$, by weak approximation of the split torus $T'_L$, it suffices to choose an $x\in T'(L)$ approximating $a$ at $w$ and approximating $1$ at the other places.
 We conclude that $b$ is approximated by $N_{L/F}(x)\in H(F)$, and hence lies in $\ov{H(F)}\cap H(F_v)\subset \prod_{v\in \cV}H(F_v)$.

The image $U$ of the norm map $N_{L_w/F_v}: T(L_w)\ra T(F_v)$ satisfies $U'=g_0Ug_0^{-1}$ for a $g_0\in V\subset H(F_v)$ such that $T\pr=g_0Tg_0^{-1}$, so each $u\in U$ is $g_0^{-1} u_0 g_0$, where $u_0\in U'$ is approximated by elements in $\ov{H(F)}\cap H(F_v)$. When shrinking $V$, we find that the associated $g_0$ is approximating the neutral element of $H(F_v)$ such that $u$ is approximated by $u_0=g_0ug_0^{-1}$. Therefore, there are elements in $\ov{H(F)}\cap H(F_v)$ approximating $u$ and $U\subset \ov{H(F)}\cap H(F_v)$.
\epf

\bprop \label{opennormalsubgroup}
There is a normal open subgroup $N$ of $\prod_{v}H(F_v)$ contained in the closure of $H(F)$.
\eprop
\bpf
It suffices to construct a normal open subgroup of $H(F_v)$ contained in $\ov{H(F)}\cap H(F_v)$ for each $v$. For a fixed maximal torus $T$ of $H_{F_v}$, by \Cref{maxtorusopensubset}, there is an open subgroup $U\ce (N_{L_w/F_v}T)(L_w)\subset T(F_v)$ contained in $\ov{H(F)}\cap H(F_v)$. Consider the following morphism
\[
f\colon H_{F_v}\times T \ra H_{F_v} \qqq (g,t)\mapsto gtg^{-1}.
\]
By the construction in the proof of \cite{SGA3II}*{XIII, 3.1}, there is a principal open dense subscheme $W\subset T$ such that every $t_0\in W(F_v)$ is a regular element, that is, by \cite{SGA3II}*{XIII, 3.0, 2.1, 2.2}, the geometric point $\ov{t_0}\in T_{\ov{F_v}}(\ov{F_v})$ over $t_0$ satisfies the following property:
\[
\text{$f_{\ov{t_0}}: H_{\ov{F_v}}\times T_{\ov{F_v}}\ra H_{\ov{F_v}} \qqq (g,t)\mapsto gtg^{-1}$ \qq is smooth at $(\id, \ov{t_0})$.}
\]
Note that this morphism is not a morphism of group schemes in general.
We claim that $U\cap W(F_v)\neq \emptyset$.
For the splitting field $L_w$ of $T$, we consider the preimage $Nr^{-1}(W)$ for the norm map $Nr\colon \mathrm{Res}_{L_{w}/F_v}(T_{L_w})\ra T$.
 Because $T_{L_w}$ is isomorphic to a dense open of an affine space $\bA_{L_w}^n$, the Weil restriction $\mathrm{Res}_{L_w/F_v}(T_{L_w})$ is also isomorphic to a dense open of $\bA_{F_v}^{mn}$ for $m=[L_w: F_v]$.
 Since $F_v$ is infinite, the Zariski density of $(\mathrm{Res}_{L_w/F_v}(T_{L_w}))(F_v)$ implies that $(\mathrm{Res}_{L_w/F_v}(T_{L_w}))(F_v)\cap Nr^{-1} (W)(F_v)\neq \emptyset$ and 
 \[
 U\cap W(F_v)=Nr((\mathrm{Res}_{L_w/F_v}(T_{L_w}))(F_v))\cap W(F_v)\neq \emptyset.
 \]
Therefore, we may choose a $t_{0}\in U\cap W(F_v)$ such that $f$ is smooth at $(\id, t_0)$. 
Hence, there is an open neighborhood $W\pr$ of $(\id, t_{0})$ such that $f|_{W\pr}$ is a smooth morphism. 
Since a smooth morphism is a composite of \'etale morphism and a projection (see \SP{054L}),
by \cite{Ces15d}*{2.8}, the map $f|_{W\pr}(F_v)\colon W\pr(F_v)\ra H(F_v)$ is open for the $v$-adic topology.
By shrinking $W\pr$ if necessary and by translation, we find an open neighborhood $U_{0}\subset U$ of $t_{0}$ such that $E\ce f(H(F_v)\times U_0)$ is open. 
Let $N$ be the subgroup of $H(F_v)$ generated by $E$. 
Since $N$ contains an open subset $E$, it is an open subgroup.
Because $E$ is stable under $H(F_v)$-conjugation, $N$ is a normal subgroup of $H(F_v)$.
By \Cref{maxtorusopensubset}, the conjugation $gtg^{-1}$ is in $\wt{U}$ for another maximal torus $\wt{T}$ of $H_{F_v}$, and hence $gtg^{-1}\in \ov{H(F)}\cap H(F_v)$. Consequently, $N$ is the desired normal open subgroup of $H(F_v)$ contained in $\ov{H(F)}\cap H(F_v)$.
\epf

\addtocontents{toc}{\protect\setcounter{tocdepth}{1}}

\section{Decomposition of groups} \label{decomposition of groups}

The goal of this section is to prove the following decomposition of groups, which leads to the reduction of \Cref{GSconj} to the case of complete discrete valuation rings (see \Cref{reducetocompletecase}):
\bthm \label{grpdecomp}
For a reductive group scheme $G$ over a semilocal Dedekind domain $R$ with fraction field $K$, let $\cV$ be the set of valuations corresponding to the maximal ideals of $R$. For each $v\in \cV$, let $R_{v}$ be the completion of $R$ at $v$ and its fraction field be $K_{v}$. Then, we have
\[
\tst \prod_{v\in \cV}G(K_v)=G(K)\prod_{v\in \cV}G(R_v).
\]
\ethm
The proof proceeds in Propositions \ref{max split tori}, \ref{paraunipotent}, \ref{paraboliccontained} and \ref{parabolicabsorb}, by mildly simplifying and improving Nisnevich's argument. A minimal parabolic subgroup of $G_{R_v}$ is denoted by $P_v$, and its unipotent radical and Levi subgroup are denoted by $U_v$ and $L_v$ respectively. By \cite{SGA3IIInew}*{XXVI, 6.11, 6.18}, the maximal central split torus $T_v$ of $L_v$ is a maximal split torus of $G_{R_v}$.

\bprop \label{max split tori}
We have
\[
\tst \prod_{v\in \cV}T_v(K_v)\subset G(K)\prod_{v\in \cV}G(R_v).
\]
\eprop
\bpf
It suffices to replace each $T_v$ by a maximal torus of $G_{R_v}$ containing it to prove the inclusion. For a fixed $v$, the morphism
\[
  f\colon G(K_v)  \ra \underline{\mathrm{Tor}}(G)(K_v) \qqq g  \mapsto gT_vg^{-1}
\]
sends a neighbourhood $U$ of the neutral element of $G(K_v)$ to {an open subset $V\subset \mtg(G)(K_v)$} containing $T_v$. Therefore, $V\cap \mtg(G)(R_v)$ is an open subset of $\mtg(G)(K_v)$ containing $T_v$. Since $\underline{\mathrm{Tor}}(G)$ is an affine $R$-scheme (\cite{SGA3II}*{XII, 5.4}), we have a Cartesian square with inclusion arrows
\[
\text{$\begin{matrix}
   \mtg(G)(R) & \ra & \mtg(G)(R_v) \\
  \downarrow &  & \downarrow \\
  \mtg(G)(K) & \ra & \mtg(G)(K_v)
\end{matrix}$\q
and \q
$\mtg(G)(R)=\mtg(G)(K)\cap \mtg(G)(R_v)$.}
\]

As in the proof of \Cref{maxtorusopensubset}, $\mtg(G)(K)$ is dense in $\mtg(G)(K_v)$, so the intersection $V\cap \mtg(G)(R_v) \cap \mtg(G)(K)=V\cap \mtg(G)(R)\neq \emptyset$ contains a maximal $R$-torus $T_0$ of $G$. By \Cref{opennormalsubgroup}, there is an open normal subgroup $N$ of $\prod_{v}G(K_v)$ contained in $\ov{G(K)}$. Assuming that $U\subset N\cap \prod_{v}G(R_v)$, we have
\[
T_v(K_v)=gT_0(K_v)g^{-1}=gT_0(K)T_0(R_v)g^{-1}\subset gG(K)G(R_v)g^{-1}\subset \ov{G(K)}G(R_v),
\]
where the second equality is by \cite{CTS87}*{Prop.~8.1}. By \cite{Con12}*{Prop.~2.1} and that $G$ is affine,  $\prod_{v}G(R_v)$ is both closed and open in $\prod_{v}G(K_v)$, so the product $G(K)\prod_{v}G(R_v)$ contains $\ov{G(K)}$. Therefore,
\[
\text{$\ov{G(K)}\prod_{v}G(R_v)\subset G(K)\prod_{v}G(R_v)$\qq and \qq $\prod_v T_v(K_v)\subset G(K)\prod_v G(R_v)$.} \qedhere
\]
\epf

\bprop \label{paraunipotent}
We have
\[
\tst \prod_{v\in \cV}U_v(K_v)\subset \ov{G(K)}.
\]
\eprop
\bpf
The maximal split torus $T_{v}$ with character lattice $M$ acts on $G_{R_v}$ by conjugation
\[
	T_v\times G_{R_v}\ra G_{R_v},\qq (t, g)\mapsto tgt^{-1},
\]
yielding a weight decomposition $\mathrm{Lie}(G_{R_v})=\bigoplus_{\alpha\in M}\Lie(G_{R_v})^{\alpha}$. 
The subset $\Phi$ of $\alpha \in M - \{0\}$ such that $\mathrm{Lie}(G_{R_v})^{\alpha}\neq 0$ is the relative root system of $(G_{R_v}, T_v)$. 
By \cite{SGA3IIInew}*{XXVI, 6.1;  7.4}, the zero-weight space of $\Lie(G_{R_v})$ is $\mathrm{Lie}(L_v)$. The relative root datum $((G_{R_v}, T_v), M, \GF)$ has a set $\GF_{+}$ of positive roots such that	
\[
\tst	\text{$\mathrm{Lie}(P_v)=\mathrm{Lie}(L_v)\oplus \Lie(U_v)=\mathrm{Lie}(L_v)\oplus \p{\bigoplus_{\alpha\in \GF_+}\mathrm{Lie}(G_{R_v})^{\alpha}}$.}
\]
Let $\wt{K_v}/K_{v}$ be a Galois field extension splitting $G_{K_v}$.
We denote the base change of $P_{v}, L_{v}, U_{v}$ and $T_{v}$ over $\wt{K_v}$ by $\wt{P}, \wt{L}, \wt{U}$ and $\wt{T}$ respectively.
Since $T_{v}$ is central in $L_{v}$, by \cite{SGA3IIInew}*{XXVI, 2.4; XIX, 1.6.2(iii)}, there is a split maximal torus $T\pr\subset \wt{P}$ of $G_{\wt{K_v}}$ containing $\wt{T}$. 
The centralizer of $T\pr$ in $G_{\wt{K_v}}$ is itself and is a Levi subgroup of a Borel subgroup $\wt{B}\subset \wt{P}$.
Let $M\pr$ be the character lattice of $T\pr$.
Then, the adjoint action of $T\pr$ on $G_{\wt{K_v}}$ induces a decomposition $\Lie(G_{\wt{K_v}})=\bigoplus_{\alpha\in M\pr}\Lie(G_{\wt{K_v}})^{\alpha}$, whose coarsening is $\Lie(G_{R_v})=\bigoplus_{\alpha\in M}\Lie(G_{R_v})^{\alpha}$. 
Let $\GF\pr$ be the root system with positive set $\GF\pr_{+}$ for the Borel $\wt{B}$. 
Now \cite{SGA3IIInew}*{XXVI, 7.12} gives us a surjection $u\colon M\pr\surjects M$ such that $\GF_{+}\subset u(\GF\pr_+)\subset \GF_{+}\cup \{0\}$.
By \emph{op.~cit.} 1.12, we have decompositions 
\[
\tst	\wt{U}=\prod_{\alpha\in \GF\prpr}\wt{U}_{\alpha},\qq \Lie(\wt{U})=\bigoplus_{\alpha\in \GF\prpr}\Lie(G_{\wt{K_v}})^{\alpha},
	\] 
where $\GF\prpr\subset \GF\pr_{+}$.
Further, we have isomorphisms of group schemes $f_{\alpha}\colon \bG_{a,\wt{K_v}}\isomto \wt{U}_{\alpha}$.
The zero-weight space for $T_{v}$-action on $\Lie(G_{R_v})$ is $\Lie(L_v)$, so the restriction to $T$ of weights in $\Lie(\wt{U})$ must be nonzero, that is $u(\GF\prpr)\subset \GF_{+}$. 
For a cocharacter $\xi\colon \bG_m\ra T_{v}$ and its base change $\xi_{\wt{K_v}}\colon \bG_{m,\wt{K_v}}\ra \wt{T}$, the composite of $\xi_{\wt{K_v}}$ with $\wt{T}\ra T\pr$ is the image $u^{\ast}(\xi)$ of $\xi$ under the dual map $u^{\ast}\colon M^{\ast}\injects {M\pr}^{\ast}$.
The action of $\bG_{m}$ on $U_{v}$ induced by $\xi$ is 
\[
\tst \mathrm{ad}\colon	\bG_m(K_v)\times U_v(K_v)\ra U_v(K_v),\qq (t,y)\mapsto \xi(t)y\xi(t)^{-1}.
\]
For a $y\in U_{v}(K_v)$, let $\wt{y}\in U_{v}(\wt{K_v})$ denote the image of $y$.
Let $N$ be the open normal subgroup of $\prod_{v}G(K_v)$ constructed in \Cref{opennormalsubgroup}.
We consider the commutative diagram
\[
\tikz {
\node  (C) at (-5.1,0) {$\bG_m(K_v)\times (N\cap U_v(K_v))$};
\node  (A) at (0,0) {$T_v(K_v)\times (N\cap U_v(K_v))$};
\node  (B) at (4.5,0) {$N\cap U_v(K_v)$};
\node  (c) at (-5.1,-1.5) {$\bG_m(K_v)\times U_v(K_v)$};
\node  (a) at (0,-1.5) {$T_v(K_v)\times U_v(K_v)$};
\node  (b) at (4.5,-1.5)   {$U_v(K_v)$};
\node  (d) at (-5.1,-3) {$\bG_m(\wt{K_v})\times U_v(\wt{K_v})$};
\node  (e) at (0,-3) {$T_v(\wt{K_v})\times U_v(\wt{K_v})$};
\node  (f) at (4.5,-3)   {$U_v(\wt{K_v})$.};
\node  (l) at (-2.5,0.2) {\small$\xi\times \id$};
\node  (ll) at (-2.5,-1.3) {\small$\xi\times \id$};
\node  (lll) at (-2.5,-2.8) {\small$\xi\times \id$};
\node  (j) at (2.6,0.2) {\small$\mathrm{ad}$};
\node  (jj) at (2.6,-1.3) {\small$\mathrm{ad}$};
\node  (jjj) at (2.6,-2.8) {\small$\mathrm{ad}$};
\draw [draw = black, thin,
arrows={
- Stealth }]
(e) edge  (f);
\draw [draw = black, thin,
arrows={
- Stealth }]
(d) edge  (e);
\draw [draw = black, thin,
arrows={
- Stealth }]
(b) edge  (f);
\draw [draw = black, thin,
arrows={
- Stealth }]
(a) edge  (e);
\draw [draw = black, thin,
arrows={
- Stealth }]
(c) edge  (d);
\draw [draw = black, thin,
arrows={
- Stealth }]
(A) edge  (B);
\draw [draw = black, thin,
arrows={
- Stealth }]
(a) edge (b);
\draw [draw = black, thin,
arrows={
- Stealth }]
(B) edge (b);
\draw [draw = black, thin,
arrows={
- Stealth }]
(A) edge (a);
\draw [draw = black, thin,
arrows={
- Stealth }]
(C) edge  (c);
\draw [draw = black, thin,
arrows={
- Stealth }]
(C) edge  (A);
\draw [draw = black, thin,
arrows={
- Stealth }]
(c) edge  (a);
}
\]
Let $\varpi$ be a uniformizer of $K_{v}^{\times}=\bG_m(K_v)$.
For an integer $n$, the action of $\varpi^{n}$ on $y$ via $\xi$ is denoted by $(\varpi^n)\cdot y$.
Decompose $\wt{y}=\prod_{\alpha\in \GF\prpr}f_{\alpha}(z_{\alpha})$ for $z_{\alpha}\in \wt{K_v}$.
The image $(u^{\ast}(\xi)(\varpi^n))\wt{y}(u^{\ast}(\xi)(\varpi^n))^{-1}\in U_{v}(\wt{K_v})$ of $(\varpi^n)\cdot y$ is the following
\[
\tst	\prod_{\alpha\in \GF\prpr}(u^{\ast}(\xi)(\varpi^n))f_{\alpha}(z_{\alpha})(u^{\ast}(\xi)(\varpi^n))^{-1}=\prod_{\alpha\in \GF\prpr}f_{\alpha}((\varpi^n)^{\langle u^{\ast}(\xi), \alpha \rangle}z_{\alpha})=\prod_{\alpha\in \GF\prpr}f_{\alpha}((\varpi^n)^{\langle \xi, u(\alpha)\rangle}z_{\alpha}).
\]
Because $u(\GF\prpr)\subset \GF_{+}$, we can choose a cocharacter $\xi$ such that $\langle \xi, u(\alpha)\rangle$ is positive for all roots $\alpha\in \GF\prpr$. 
Then, when $n$ grows, the element $(\varpi^n)\cdot y\in U_{v}(\wt{K_v})$ gets closed to identity and so the same holds in $U_{v}(K_v)$.
Thus, since $N\cap U_{v}(K_v)$ is an open neighbourhood of $\mathrm{id}_{G(K_v)}$, every orbit of the $T_{v}$-action on $U_{v}(K_v)$ intersects with $N\cap U_{v}(K_v)$ nontrivially, so
\[
\tst	U_v(K_v)=\bigcup_{t\in T_v(K_v)}t(N\cap U_v(K_v))t^{-1}=N\cap U_v(K_v),
\]
which implies that $U_{v}(K_v)\subset N$.
Combining \Cref{opennormalsubgroup}, we conclude that $\prod_{v\in \cV}U_v(K_v)\subset \ov{G(K)}$.\qedhere
\epf
Before the next step of decomposing the group, we gather
some recollections on anisotropic groups.
\blem\label{genericanisotropic}
   For a reductive group $G$ over a discrete valuation ring $R$ with fraction field $K$,
	\[
	\text{$G$ is $R$-anisotropic \q if and only if \q $G$ is $K$-anisotropic.}
	\]
\elem
\bpf
For an $R$-algebra $R\pr$, by \cite{SGA3IIInew}*{XXVI, 6.14}, the group $G$ is $R'$-anisotropic if and only if $G_{R'}$ has no proper parabolic subgroups and $\rad(G_{R'})$ has no split subtorus. The proof proceeds in the following two steps.
  \benum
	\item[(a)] For the existence of parabolic subgroups, by \cite{SGA3IIInew}*{XXVI, 3.5}, the functor
	 \[\qq
	 \underline{\mathrm{Par}}(G)\colon  \mathbf{Sch}_{/R} \longrightarrow \mathbf{Sets}, \qqq
	   R'  \longmapsto  \{\text{parabolic subgroups of $G_{R'}$}\}
	 \]
is representable by a projective scheme over $R$ and satisfies the valuative criterion of properness:
  \[
  \qq \underline{\mathrm{Par}}(G)(K)=\un{\mathrm{Par}}(G)(R).
  \]
	\item[(b)] If $\rad(G_K)=\rad(G)_K$ is $K$-anisotropic, then $\rad(G)$ is $R$-anisotropic. 
	For the other direction, recall that every torus over a Noetherian normal local ring $S$ is isotrivial, see \cite{SGA3II}*{X, 5.16}.
    We use the anti-equivalence in \cite{SGA3II}*{X, 1.2}:
\[\qq \hkh{\begin{matrix}
			   \text{The category of}  \\
			   \text{tori over $S$}
			 \end{matrix} } \longleftrightarrow
	  \hkh{\begin{matrix}
			 \text{The category of $\mathbb{Z}$-lattices of finite type} \\
			 \text{with continuous action of the fundamental group of $S$}
		   \end{matrix} }.
\]

Both $\rad(G)$ and $\rad(G_K)$ correspond to a $\mathbb{Z}$-lattice $M$ acted by \'etale fundamental groups $\pi^{\et}_{1}(R)$ and $\pi_{1}^{\et}(K)$ respectively. A nontrivial split torus in $\rad(G_K)$ corresponds to a quotient lattice $N$ of $M$ with trivial $\pi_1^{\et}(K)$-action. By the surjectivity of $\pi_1^{\et}(K)\to \pi_{1}^{\et}(R)$(see \cite[\href{https://stacks.math.columbia.edu/tag/0BSM}{0BSM}]{SP}), the action of $\pi_{1}^{\et}(R)$ on $N$ is also trivial and the latter corresponds to a nontrivial split torus in $\rad(G)$. \qedhere
\eenum
\epf
\bprop \label{anisotropicrationalintegerpoints}
  For a reductive group $G$ over a Henselian discrete valuation ring $R$,
  \benum
 \item for the strict Henselization $\wt{R}$ of $R$ with fraction field $\wt{K}$, the group $G(\wt{R})$ is a maximal parahoric subgroup of $G(\wt{K})$\uscolon
  \item if $G$ is $K$-anisotropic, then $G(K)=G(R)$.
  \eenum
\eprop
\bpf \hfill
\benum
\item

Because $G_{\wt{R}}$ is a Chevalley group scheme and the residue field of $\wt{R}$ is separably closed, by \cite{BrT2}*{4.6.22, 4.6.31}, the group $G(\wt{R})$ is the stabilizer of a special point in the Bruhat--Tits building of $G(\wt{K})$. Since a special point is a minimal facet (see \cite{BrT2}*{4.6.15}), by connectedness of each fiber of $G_{\wt{R}}/\wt{R}$ and the definition of parahoric subgroups (see \cite{BrT2}*{5.2.6}), $G(\wt{R})$ is a maximal parahoric subgroup of $G(\wt{K})$.
\item First, we prove that $G(K)\subset G(\wt{R})$. For the Bruhat--Tits building $\wt{\sI}$ of $G(\wt{K})$, by \cite{BrT2}*{5.1.27}, the Galois group $\Sigma\ce \Gal(\wt{K}/K)$ acts on the enlarged Bruhat--Tits building $\wt{\sI}^{ext}\ce \wt{\sI}\times X^{\ast}_{\wt{K}}(G)_{\b{R}}^{\vee}$ of $G(\wt{K})$ and has a unique fixed point $x$ located at $\wt{\sI}\times \{0\}$. For $x$, by \cite{BrT2}*{5.2.6}, its connected pointwise stabilizer (see \cite{BrT2}*{4.6.28}) is a $\Sigma$-invariant parahoric subgroup $P$. This parahoric subgroup $P$, by the $K$-anisotropicity of $G$ and \cite{BrT2}*{5.2.7}, is the unique $\Sigma$-invariant parahoric subgroup and thus, due to (a), coincides with $G(\wt{R})$. By connectedness of each fiber of $G_{\wt{R}}/\wt{R}$, the parahoric subgroup $G(\wt{R})$ is also the stabilizer of $x$. By \cite{BrT1}*{9.2.1 (DI 1)}, the fixed point $x$ is also stabilized by $G(K)$. Therefore, $G(K)\subset G(\wt{R})$.

	Lastly, by the affineness of $G$ and $R=\wt{R}\cap K$, we obtain the following cartesian square
	\[
	\begin{matrix}
   G(R) & \ra & G(\wt{R}) \\
  \downarrow &  & \downarrow \\
   G(K) & \ra & G(\wt{K})
\end{matrix}
	\]
	 and combine it with $G(K)\subset G(\wt{R})$ to conclude that $G(K)=G(R)$. \qedhere
\eenum
\epf
\brem
When using Bruhat--Tits theory for general reductive groups, one needs to check the conditions listed in \cite{BrT2}*{5.1.1}. Fortunately, these conditions are automatically satisfied when $R$ is a Henselian discrete valuation ring, and $G$ is defined and reductive over $R$ (by \cite{SGA3IIInew}*{XXII, 2.3}, $G$ splits over $\wt{R}$, hence is split (a fortiori quasi-split) over $\wt{K}$).
\erem
\bprop \label{paraboliccontained}
We have
\[
\tst \prod_{v\in \cV}P_v(K_v)\subset G(K)\prod_{v\in \cV}G(R_v).
\]
\eprop
\bpf
The quotient $H_v\ce L_v/T_v$ is $R_v$-anisotropic and by \Cref{genericanisotropic} is $K_v$-anisotropic. By \Cref{anisotropicrationalintegerpoints}, we have $H_v(K_v)=H_v(R_v)$, which fits into the commutative diagram with exact rows:
\[
\tikz {
\node  (A) at (0,0) {$0$};
\node  (B) at (2,0) {$T_v(R_v)$};
\node  (C) at (4,0) {$L_v(R_v)$};
\node  (D) at (6,0) {$H_v(R_v)$};
\node  (E) at (8.5,0) {$H^1(R_v, T_v)=0$};
\node  (F) at (0,-1.5)   {$0$};
\node  (H) at (2,-1.5) {$T_v(K_v)$};
\node  (I) at (4,-1.5) {$L_v(K_v)$};
\node  (J) at (6,-1.5) {$H_v(K_v)$};
\node  (K) at (8.5,-1.5) {$H^1(K_v, T_v)=0$};
\node  (L) at (4.25,-0.73) {$\lambda_v$};
\draw [draw = black, thin,
arrows={
- Stealth }]
(A) edge  (B);
\draw [draw = black, thin,
arrows={
- Stealth }]
(B) edge  (C);
\draw [draw = black, thin,
arrows={
- Stealth }]
(C) edge (D);
\draw [draw = black, thin,
arrows={
- Stealth }]
(D) edge (E);
\draw [draw = black, thin,
arrows={
- Stealth }]
(F) edge (H);
\draw [draw = black, thin,
arrows={
- Stealth }]
(H) edge (I);
\draw [draw = black, thin,
arrows={
- Stealth }]
(I) edge (J);
\draw [draw = black, thin,
arrows={
- Stealth }]
(J) edge (K);
\draw [draw = black, thin,
arrows={
- Stealth }]
(B) edge (H);
\draw [draw = black, thin,
arrows={
- Stealth }]
(C) edge (I);
\draw (5.96cm,-3.2em) -- (5.96cm,-0.8em);
\draw (6.04cm,-3.2em) -- (6.04cm,-0.8em);
\draw [draw = black, thin,
arrows={
- Stealth }]
(E) edge (K);
}.
\]
The image of $a\in L_v(K_v)$ in $H_v(K_v)$ has a preimage $b\in L_v(R_v)$. Hence $a\cdot \lambda_v(b)^{-1}\in T_v(K_v)$ and
\[
L_v(K_v)=T_v(K_v)L_v(R_v).
\]
By \Cref{max split tori} and \Cref{paraunipotent}, we obtain the conclusion.
\epf

\bprop \label{parabolicabsorb}
  \[
  \tst \ov{G(K)}\prod_{v\in \cV}P_v(K_v)=\prod_{v\in \cV}G(K_v).
  \]
\eprop
\bpf
For each $P_v$, there is a unique parabolic subgroup $Q_v$ of $G_{R_v}$ such that the canonical morphism
	\[
	\rad^u(P_v)(K_v)\cdot \rad^u(Q_v)(K_v)\to G(K_v)/P_v(K_v)
	\]
	is surjective (\cite{SGA3IIInew}*{XXVI, 4.3.2, 5.2}). We conclude by combining this surjectivity with \Cref{paraunipotent} and deducing
\[
\tst \prod_{v\in \cV}G(K_v)\subset \ov{G(K)}\prod_{v\in \cV}P_v(K_v). \qedhere
\]
\epf
\bpf[The proof of \Cref{grpdecomp} \upshape:]
By \Cref{paraboliccontained} and \Cref{parabolicabsorb}, we have
\[
  \tst    \prod_{v\in \cV}G(K_v)\subset \ov{G(K)}\prod_{v\in \cV}G(R_v)=G(K)\prod_{v\in \cV}G(R_v). \qedhere
	  \]
\epf

\addtocontents{toc}{\protect\setcounter{tocdepth}{1}}

\section{Reductions: twisting, passage to completion, and to anisotropic groups}\label{reduction of the GS conj}
This section exhibits a sequence of methods to reduce \Cref{GSconj}. First, using twisting technique of torsors, one can prove \Cref{GSconj} by showing that the map
\[
H^1_{\et}(R, G)\to H^1_{\et}(K, G)
\]
has trivial kernel. Secondly, we improve and extend Nisnevich's argument by a patching technique for the semilocal case and reduce to the case when $R$ is a complete discrete valuation ring (\Cref{reducetocompletecase}). One motivation for reducing to the complete case is that in Bruhat--Tits theory \cite{BrT2}*{\S5}, to define parahoric subgroups for more general cases than the quasi-split case, we need additional conditions (see \cite{BrT2}*{5.1.6}), which are automatically satisfied in our case when $R$ is complete. Lastly, with crucial results and properties derived from Bruhat--Tits theory (for instance, \Cref{anisotropicrationalintegerpoints}), we reduce further to the case when $G$ is semisimple and anisotropic (\Cref{reducetosemisimpleanisotropic}). This case has enough advantages such as the uniqueness of Galois-invariant parahoric subgroups of $G(\wt{K})$ (see \cite{BrT2}*{5.2.7}), for us to prove \Cref{GSconj} in \S\ref{proofssanisotropic}.

\bprop[\cite{Gir71}*{III, 2.6.1~(i)}]\label{kestutis}
For a group scheme $G$ over a scheme $S$, we let $T$ be a right $G$-torsor and let $_TG\ce \mathrm{Aut}_G(T)$. Then twisting by $T$ induces an isomorphism
\[
	  H^1_{\et}(S, G)  \stackrel{\sim}{\longrightarrow} H^1_{\et}(S, \leftsub{T}{G} ), \qqq
	   X  \longmapsto   \mathrm{Hom}_G(T,X)
\]
such that the image of the class of $T$ is the neutral class. In fact, there is an equivalence between the gerbe of $G$-torsors and the gerbe of $\leftsub{T}{G}$-torsors.
\eprop


\begin{cor} \label{GSconjtrivialkernel}
  \Cref{GSconjgeneral} is equivalent to showing that for all reductive groups $G$ over $R$,
  \[
  \text{$H^1_{\et}(R, G)\to H^1_{\et}(K, G)$ \q is injective.}
   \]
In particular, in \Cref{GSconj}, it suffices to prove that the map $H^1_{\et}(R,G)\ra H^1_{\et}(K,G)$ has trivial kernel.
\end{cor}
The following shows that we can reduce \Cref{GSconj} to the complete case.
\bprop \label{reducetocompletecase}
It suffices to prove \Cref{GSconj} when $R$ is a complete discrete valuation ring.
\eprop
\bpf
	If a $G$-torsor $\cX$ becoming trivial over $K$, then it becomes trivial over $\coprod_{v} \Spec K_v$ and by assumption, trivial over $\coprod_{v} \Spec R_v$. Because $\coprod_{v}\Spec R_v\ra \Spec R$ is flat and induces isomorphisms on closed points, by patching technique in \cite{MB14}*{Thm.~1.1}, we have an equivalence of categories
	\[\q \hkh{\begin{matrix}
			   \text{affine schemes}  \\
			   \text{over $\Spec R$}
			 \end{matrix} } \longleftrightarrow
	  \hkh{\begin{matrix}
			 \text{$(\cX\pr,\cX\prpr,\iota)$ with  $\cX\pr$ an affine scheme over $\coprod_{v}\Spec R_{v}$, $\cX\prpr$ an affine} \\
			 \text{scheme over $\Spec K$ and $\iota: \cX\pr|_{\coprod_{v}\Spec K_{v}}\isomto \cX\prpr|_{\coprod_{v}\Spec K_{v}}$}
		   \end{matrix} }.
\]
Let $\cX\pr$ be a $G_{\coprod_v\Spec R_v}$-torsor and $\cX\prpr$ a $G_{\Spec K}$-torsor.
Since $G$ is affine and flat, the resulting scheme $\wt{\cX}$ is a $G$-torsor with the glued structural isomorphism
\[
G\times_R \wt{\cX} \isomto \wt{\cX}\times_R \wt{\cX},\qq (g,s)\mapsto (gs,s).	
\]
Now let $\cT$ be an arbitrary $G$-torsor obtained by  gluing the trivial torsors ${\cX_{\coprod _{v} \Spec R_v}}$ and ${\cX_{\Spec K}}$  along ${\cX_{\coprod_{v} \Spec K_v}}$.
The torsor class of $\cT$ is determined by the gluing isomorphism $\iota: (\cX_{\coprod _{v} \Spec R_v})|_{\coprod_v\Spec K_v}\isomto (\cX_{\Spec K})|_{\coprod_v\Spec K_v} $.
This isomorphism between trivial torsors is nothing but a multiplication by $a\in \prod_vG(K_v)$.
Since automorphisms of $\cX_{\Spec K}$ and of $\cX_{\coprod_v\Spec R_v}$ are induced by multiplications by elements in $G(K)$ and in $\prod_vG(R_v)$ respectively, the isomorphic classes of the possible $G$-torsors $\cT$ are in correspondence with
\[
	\textstyle G(K)\backslash \prod_{v}G(K_v) / \prod_{v}G(R_v),
\]
which is trivial by \Cref{grpdecomp}.
Consequently, $\cX$ is in the class of trivial torsors.    
\epf
Now we reduce \Cref{GSconj} to the case when $G$ is anisotropic and semisimple.
\bprop \label{reducetonoproperparabolic}
Let $R$ be a valuation ring with fraction field $K$ and let $G$ be a reductive $R$-group scheme.
If for every reductive $R$-group scheme $H$ without a proper parabolic subgroup, the map $H^{1}_{\et}(R,H)\ra H^{1}_{\et}(K,H)$ is injective, then so is $H^{1}_{\et}(R,G)\ra H^{1}_{\et}(K,G)$. 
\eprop
\bpf
If there is a proper minimal parabolic subgroup $P$ of $G$ and $L$ is its Levi subgroup, then we consider the following commutative diagram
\[
\tikz {
\node  (C) at (-3.2,0) {$H^1_{\et}(R,L)$};
\node  (A) at (-0.6,0) {$H^1_{\et}(R,P)$};
\node  (B) at (2,0) {$H^1_{\et}(R,G)$};
\node  (c) at (-3.2,-1.5) {$H^1_{\et}(K,L)$};
\node  (a) at (-0.6,-1.5) {$H^1_{\et}(K,P)$};
\node  (b) at (2,-1.5)   {$H^1_{\et}(K,G)$};
\draw [draw = black, thin,
arrows={
- Stealth }]
(A) edge  (B);
\draw [draw = black, thin,
arrows={
- Stealth }]
(a) edge (b);
\draw [draw = black, thin,
arrows={
- Stealth }]
(B) edge (b);
\draw [draw = black, thin,
arrows={
- Stealth }]
(A) edge (a);
\draw [draw = black, thin,
arrows={
- Stealth }]
(C) edge  (c);
\draw [draw = black, thin,
arrows={
- Stealth }]
(C) edge  (A);
\draw [draw = black, thin,
arrows={
- Stealth }]
(c) edge  (a);
}
\]
and show that the kernel of the third column comes from the kernel of the first column as follows.

By \cite{SGA3IIInew}*{XXVI, 2.3}, the rows in the left square are isomorphisms.
For the second square, we use the argument in \cite{nisnevich1983etale}*{Prop.~5.1}. Let $E$ be a $G$-torsor that becomes trivial over $K$.
The quotient sheaf $E/P$ is fpqc locally isomorphic to $G/P$ and represented by a scheme projective over $R$ (see \cite{SGA3IIInew}*{XXVI, 3.3; 3.20}).
By the valuative criterion of properness, we have
\[
(E/P)(K)=(E/P)(R).
\]
For a $K$-point $x$ of $E$, its image $\ov{x}$ in $E/P$ is also an $R$-point of $E/P$. Subsequently, the fiber of $E$ over $\overline{x}$
\[
F\ce E\times_{E/P}\Spec(R)
\]
is a $P$-torsor over $R$ such that $F(K)\neq \emptyset$ and its image under the map $H^1_{\et}(R,P)\to H^1_{\et}(R,G)$ is the class of $E$. Therefore, the kernel of $H^1_{\et}(R,G)\to H^1_{\et}(K,G)$ is in the image of the kernel of $H^1_{\et}(R,P)\to H^1_{\et}(K,P)$.

By \cite{SGA3IIInew}*{XXVI, 1.20}, parabolic subgroups of $L$ are intersections of $L$ with parabolic subgroups contained in $P$. Therefore, $L$ contains no proper parabolic subgroup.
\epf

\bprop\label{reducetosemisimpleanisotropic}
In order to prove \Cref{GSconj}, it suffices to prove that $H^{1}_{\et}(R,G)\ra H^{1}_{\et}(K,G)$ has trivial kernel when $R$ is a complete discrete valuation ring and $G$ is semisimple and $R$-anisotropic.
\eprop

\bpf
By \Cref{reducetocompletecase} and \Cref{reducetonoproperparabolic}, it suffices to prove \Cref{GSconj} in the case when $G$ is a reductive group over a complete discrete valuation ring $R$ and has no nontrivial parabolic subgroup. The quotient $G/\rad(G)$ has no proper parabolic subgroups and is semisimple, hence by \cite{SGA3IIInew}*{XXVI, 6.14} is $R$-anisotropic and by \Cref{genericanisotropic} is $K$-anisotropic. We assume that $l(G/\rad(G))\colon H^1_{\et}(R, G/\rad(G))\ra H^1_{\et}(K, G/\rad(G))$ is injective. Since $R$ is complete, by \Cref{anisotropicrationalintegerpoints}, we have
\[
(G/\rad(G))(R)=(G/\rad(G))(K)
\]
fitting into the following commutative diagram
	 \[
\tikz {
\node  (A) at (-0.6,0) {$(G/\rad(G))(R)$};
\node  (B) at (3,0) {$H^1_{\et}(R,\rad(G))$};
\node  (C) at (6,0) {$H^1_{\et}(R,G)$};
\node  (D) at (9,0) {$H^1_{\et}(R,G/\rad(G))$};
\node  (a) at (-0.6,-1.5) {$(G/\rad(G))(K)$};
\node  (b) at (3,-1.5)   {$H^1_{\et}(K,\rad(G))$};
\node  (c) at (6,-1.5) {$H^1_{\et}(K,G)$};
\node  (d) at (9,-1.5) {$H^1_{\et}(K,G/\rad(G))$};
\node  (l) at (3.7,-0.75)  {$\scriptstyle{l(\rad(G))}$};
\node  (ll) at (6.4,-0.75) {$\scriptstyle{l(G)}$};
\node  (lll) at (9.9,-0.75) {$\scriptstyle{l(G/\rad(G))}$};
\draw [draw = black, thin,
arrows={
- Stealth }]
(A) edge  (B);
\draw [draw = black, thin,
arrows={
- Stealth }]
(B) edge  (C);
\draw [draw = black, thin,
arrows={
- Stealth }]
(C) edge (D);
\draw [draw = black, thin,
arrows={
- Stealth }]
(a) edge (b);
\draw [draw = black, thin,
arrows={
- Stealth }]
(b) edge (c);
\draw [draw = black, thin,
arrows={
- Stealth }]
(c) edge (d);
\draw [draw = black, thin,
arrows={
- Stealth }]
(B) edge (b);
\draw [draw = black, thin,
arrows={
- Stealth }]
(C) edge (c);
\draw (-0.64cm,-3.2em) -- (-0.64cm,-0.8em);
\draw (-0.56cm,-3.2em) -- (-0.56cm,-0.8em);
\draw [draw = black, thin,
arrows={
- Stealth }]
(D) edge (d);
}
\]
with exact rows and the map $l(\rad(G))$ is injective by \cite{CTS87}*{Thm.~4.1}. By diagram chasing, the map $l(G)$ has trivial kernel. By \Cref{GSconjtrivialkernel}, we complete the proof.
\epf

\addtocontents{toc}{\protect\setcounter{tocdepth}{2}}

\section{The semisimple and anisotropic case} \label{proofssanisotropic}
By \Cref{reducetosemisimpleanisotropic}, it suffices to prove \Cref{GSconj} for the case when $G$ is semisimple anisotropic and $R$ is a complete discrete valuation ring.
When $R$ has perfect residue field, this case is contained in \cite{nisnevich1983etale}*{2, Thm.~4.2} for general reductive groups.
However, the proof has several unclear points, see Remarks \ref{non-self-normalizing}-\ref{Nisdefectproduct}.
These gaps motivate us to first reduce to the semisimple and anisotropic case, where we have the uniqueness of Galois-invariant parahoric subgroups. 
Further, in our proof, we add supplementary details for normalizers of parahoric subgroups.
\Cref{finalproof} also extends a special case implied by \cite{BrT3}*{Lem.~3.9} when the residue field is perfect and $G$ is semisimple simply connected.

\bthm \label{finalproof}
  For a semisimple and anisotropic group scheme $G$ over a complete discrete valuation ring $R$ with fraction field $K$, the following map has trivial kernel\upshape:
  $$H^1_{\et}(R, G)\to H^1_{\et}(K, G).$$
\ethm
\bpf
We fix some notations:
\begin{itemize}
  \item let $K^{\mathrm{sep}}$ be a separable closure of $K$ that contains $\wt{K}$;
  \item
	  $\text{$\Gamma\ce \Gal(\wt{R}/R)\simeq\Gal(\wt{K}/K)$; \qqqq $\Gamma_{\wt{K}}\ce \Gal(K^{\mathrm{sep}}/\wt{K})$; \qqqq $\Gamma_K\ce \Gal(K^{\mathrm{sep}}/K)$}$.
\end{itemize}
Since $R$ is Henselian, we can view $G$ as a sheaf over the site of profinite $\GG$-sets, whose one-point set with trivial $\GG$-action is $\Spec(R/\fm_{R})$.
 Therefore, a variant of the Cartan--Leray spectral sequence \cite{Sch13}*{3.7(iii)} gives an isomorphism to Galois cohomology $H^1_{\et}(R,G)\isom H^1(\GG,G(\wt{R}))$.
Similarly (or by \cite{SGA4II}*{VIII, 2.1}), we have $H^1_{\et}(K, G)\isom H^1(\Gamma_K, G(K^{\mathrm{sep}}))$.
It suffices to prove that both $\alpha$ and $\beta$ in the decomposition
\[
H^1(\Gamma, G(\wt{R}))\stackrel{\alpha}{\to}H^1(\Gamma, G(\wt{K}))\stackrel{\beta}{\to} H^1(\Gamma_K, G(K^{\mathrm{sep}}))
\]
have trivial kernel. We show this in the following two steps.
\begin{itemize}
  \item[(i)] The injectivity of
  \[
  \alpha\colon H^1(\Gamma, G(\wt{R}))\to H^1(\Gamma, G(\wt{K})).
  \]
  For a cocycle $z\in H^1(\Gamma, G(\wt{R}))$ that becomes trivial in $H^1(\Gamma, G(\wt{K}))$, there is $h\in G(\wt{K})$ such that
  \[
  \text{$z(s)=h^{-1}s(h)$ \ in \ $G(\wt{R})$, \q for each $s\in \Gamma$.}
  \]
  To prove that $h\in G(\wt{R})$, we consider the subgroups $G(\wt{R})$ and $hG(\wt{R})h^{-1}$ of $G(\wt{K})$. We have seen that $G(\wt{R})$ is a parahoric subgroup of $G(\wt{K})$ in \Cref{anisotropicrationalintegerpoints}. The conjugation by $g\in G(\wt{K})$ of a parahoric subgroup $P_F$ associated to the facet $F$ satisfies the definition of the parahoric subgroup associated to the facet $g\cdot F$ (see \cite{BrT2}*{5.2.6}). In particular, $hG(\wt{R})h^{-1}$ is also a parahoric subgroup. Now we show that $hG(\wt{R})h^{-1}$ is invariant under $\Gamma$: for every $s\in \Gamma$,
  \[
  s(hG(\wt{R})h^{-1})=s(h)G(\wt{R})s(h^{-1})=hz(s)G(\wt{R})z(s)^{-1}h^{-1}=hG(\wt{R})h^{-1},  \]
   since $G(\wt{R})$ is $\Gamma$-invariant and $z(s)\in G(\wt{R})$. Because $G$ is anisotropic over $K$, the uniqueness of a $\Gamma$-invariant parahoric subgroup (\cite{BrT2}*{5.2.7}) implies that
\[
hG(\wt{R})h^{-1}=G(\wt{R}).
\]
By \cite{BrT2}*{4.6.22, 4.6.31}, because the residue field of $\wt{R}$ is separably closed and $G$ is a Chevalley group scheme over $\wt{R}$, the group $G(\wt{R})$ is the group of $\wt{R}$-points of a connected pointwise stabilizer $\mathfrak{G}^0_{x}$ of a special point $x$ (which is a minimal facet of the extended Bruhat--Tits building $\wt{\sI}^{ext}\ce \wt{\sI}\times X^{\ast}_{\wt{K}}(G)_{\b{R}}^{\vee}$) and is the stabilizer $\mathfrak{G}_x^{\dagger}$ of $x$ in $\wt{\sI}$ (see \cite{BrT2}*{4.6.28}). On the other hand, since $G$ is semisimple, by \cite{SGA3IIInew}*{XXII, 6.2.1}, the $\wt{K}$-character group of $G$ is $X^{\ast}_{\wt{K}}(G)\ce \mathrm{Hom}_{\text{$\wt{K}$-gr.}}(G, \mathbb{G}_{m,\wt{K}})\isom \mathrm{Hom}_{\text{$\wt{K}$-gr.}}(G/G^{der},\mathbb{G}_{m,\wt{K}})=1$. Subsequently, $\wt{\sI}^{ext}=\wt{\sI}$ and $x$ is a point in the Bruhat--Tits building, in which the stabilizer of $x$ in $G(\wt{K})$ is the normalizer of $G(\wt{R})$ in $G(\wt{K})$ (see \cite{BrT2}*{4.6.17}), so $G(\wt{R})=\mathrm{Norm}_{G(\wt{K})}(G(\wt{R}))$. Thus, we have $h\in G(\wt{R})$ and $z$ is trivial in $H^1(\Gamma, G(\wt{R}))$.

 \item[(ii)] The injectivity of
 \[
 \beta\colon H^1(\Gamma, G(\wt{K}))\to H^1(\Gamma_K, G(K^{\mathrm{sep}})).
 \]

  Recall the inflation-restriction exact sequence in \cite{Ser02}*{5.8~a)}:
  \[
  0\to H^1(G_1/G_2, A^{G_2})\to H^1(G_1, A)\to H^1(G_2,A)^{G_1/G_2},
  \]
  where $G_1$ is a group with a closed normal subgroup $G_2$ and $A$ is a $G_1$-group. It suffices to take
  \[
  \text{$G_1\ce \Gamma_K$, \q $G_2\ce \Gamma_{\wt{K}}$,\q and $A\ce G(K^{\mathrm{sep}})$.}  \qedhere
  \]
\end{itemize}
\epf
The following remarks explain why our approach to proving \Cref{GSconj} in the case of a complete discrete valuation ring would not have worked without reducing first to the semisimple anisotropic case. 
 They also highlight several problems in the argument of \cite{nisnevich1983etale}*{2, Thm.~4.2}.

\brems 
\remi \label{non-self-normalizing}
When $G$ is a semisimple, simply-connected group over a Henselian discrete valuation field $K$ (for instance, a local field), parahoric subgroups of $G(K)$ are their own normalizers (\cite{BrT2}*{5.2.9}). However, this is not true in the semisimple adjoint case. If $G\ce PGL_2(\b{Q}_p)$, then
\[
\text{\q the parahoric subgroup $P=\hkh{\text{classes of ${\begin{pmatrix}
								\b{Z}_p^{\times} & p\b{Z}_p \\
								\b{Z}_p & \b{Z}_p^{\times}
							  \end{pmatrix}}$}}$ is normalized by $A=\begin{pmatrix}
																 0 & 1 \\
																 p & 0
															   \end{pmatrix}$
}
\]
but $A\not \in P$: the diagonals of $P$ are nonzero. This also provides an example of the fact that maximal bounded subgroups are not necessarily stabilizers of vertices; the converse is true (see \cite{JKY}*{p.~14}).
\remi Maximal parahoric subgroups are not their own normalizers in the reductive case. For instance, we consider $GL_n(\b{Q}_p)$ with a maximal parahoric subgroup $GL_n(\b{Z}_p)$. The normalizer of $GL_n(\b{Z}_p)$ contains $\b{Q}_p^{\times}\cdot \mathbf{1}$, which is not in $GL_n(\b{Z}_p)$. To explain this, we note that $GL_n(\b{Z}_p)$ stabilizes $x\times \{0\}\subset \sI^{ext}=\sI\times X^{\ast}_{\b{Q}_p}(GL_n)_{\b{R}}^{\vee}$, but the stabilizer of the special point $y\ce x\times X^{\ast}_{\b{Q}_p}(GL_n)_{\b{R}}^{\vee}$ contains $\b{Q}_{p}^{\times}\cdot \mathbf{1}$, which acts on $y$ by translation only in the second component. This motivates us to reduce to the semisimple case before proving \Cref{finalproof}: the reduced Bruhat--Tits building coincides with the extended one (see \cite{BrT2}*{4.2.16}), so that $G(\wt{R})$ is its own normalizer.
\remi \label{Nisdefectproduct}
In the proof of \cite{nisnevich1983etale}*{Thm.~4.2}, the author denoted by $T$ the central torus of $G$ and by $G^{der}$ the derived group of $G$, saying that $G$ is the ``almost direct product of $T$ and $G^{der}$'', without further information about $T$. There is certainly an isogeny $G^{der}\times \rad(G)\to G$ (\cite{SGA3IIInew}*{XXII, 6.2.4}). Nevertheless, the equation
\[
G(\wt{K})=G^{der}(\wt{K})\cdot \rad(G)(\wt{K})
\]
used in \cite{nisnevich1983etale}*{Thm.~4.2} fails in general. For instance, we take $G=\mathrm{GL}_n$ and $\wt{K}=\mathbb{Q}_p^{ur}$ to find
\[
\text{$\mathrm{GL}_n(\mathbb{Q}_p^{ur})\neq \mathrm{SL}_n(\mathbb{Q}_p^{ur})\cdot \{\mathrm{diag}(a,\cdots,a)\}$, \qq where $a$ ranges over $\mathbb{Q}_p^{ur\times}$.}
\]
 In fact, for $a\in\mathbb{Q}_p^{ur\times}\backslash (\mathbb{Q}_p^{ur\times})^n$, e.g. $a=p$, the matrix $\mathrm{diag}(a,1,\cdots, 1)$ has determinant $a$ and is not a product of matrices on the right.
\erems
As an application of the Grothendieck--Serre conjecture for the case of discrete valuation rings, we prove the following proposition.
\bprop
For a reductive group scheme $G$ over a field $k$, we let $k(\!(X)\!)$ be the field of Laurent power series in the variable $X$. Then the following map is injective:
\[
H^1(k, G)\ra H^1(k(\!(X)\!), G).
\]
\eprop
\bpf
The projection $p\colon \Spec k\fps{X}\ra \Spec k$ has a section $s\colon \Spec k\ra \Spec k\fps{X}$ given by $X=0$. We have the composition $H^1(k, G)\stackrel{p^{\ast}}{\ra} H^1(k\fps{X}, G)\stackrel{s^{\ast}}{\ra} H^1(k, G)$ such that $s^{\ast}\circ p^{\ast}=\mathrm{id}$. So $p^{\ast}\colon H^1(k, G)\ra H^1(k\fps{X}, G)$ is injective (in fact, by \cite{SGA3IIInew}*{XXIV, 8.1}, we have $H^1(k,G)\isom H^1(k\llbracket X \rrbracket, G)$). Since $k\llbracket X \rrbracket $ is a discrete valuation ring, by \Cref{GSconj}, the map $H^1(k\fps{X}, G)\ra H^1(k\lps{X}, G)$ is injective.
\epf

\section{Uniqueness of reductive models}
As an application of \Cref{GSconj}, we consider a scheme $X$ with function field $K$, over which a reductive group $G$ is defined. We call a \emph{model} of $G$ a flat, affine, and finite type $X$-group scheme $\sG$ such that $\sG_K\isom G$. The question is, how many reductive models does $G$ have? When $X$ is local strictly Henselian, the uniqueness of models is a special case of \cite{PY06}*{Thm.~8.5}, where they classified quasi-reductive models. In fact, in the following, we show that if a variant of \Cref{GSconjgeneral} holds, then a reductive group scheme over a semilocal regular base is determined by its generic fiber. In particular, \Cref{modelsemilocalDedekind} implied by \Cref{Uniqueness} generalizes the special case in \cite{Nis84}*{Thm.~5.1} where $X$ is the spectrum of a discrete valuation ring and $G$ is semisimple.
\bprop\label{Uniqueness}
For a regular semilocal ring $S$ with total ring of fractions $K$, if for each reductive group scheme $G'$ over $S$, the map $H^1(S, G')\ra H^1(K, G')$ is injective, then any reductive $K$-group $G$ has at most one reductive $S$-model.
\eprop
\bpf
For a reductive group scheme $\sG$ over $S$, by \cite{SGA3IIInew}*{XXIV, 1.17}, the following functor defines a one-to-one correspondence (up to isomorphisms) of sets:
\[\hkh{\begin{matrix}
			   \text{$S$-group schemes }  \\
			   \text{that are fpqc locally isomorphic to $\sG$}
			 \end{matrix} } \rightarrow
	 H^1_{\et}(S, \un{\mathrm{Aut}}_{\text{$S$-gr.}}(\sG)) \qq \sG'\mapsto \un{\mathrm{Isom}}_{\text{$S$-gr.}}(\sG, \sG').
\]
Let $\sG$ and $\sG\pr$ be two reductive $S$-models of $G$. By \cite{SGA3IIInew}*{XXII, 2.8}, the root datum of $\sG$ at each fiber of $\sG\ra\Spec S$ is locally constant, so it is the root datum of $G=\sG_K$. Therefore, $\sG$ and $\sG\pr$ have the same root datum. By \cite{SGA3IIInew}*{XXIII, 5.1}, \'etale locally we have $\sG\isom \sG\pr$ and $\sG$ corresponds to $0\in H^1_{\et}(S, \underline{\mathrm{Aut}}_{\text{$S$-gr.}}(\sG))$. We let $\alpha$ be the class of $\sG'$ in $H^1_{\et}(S, \underline{\mathrm{Aut}}_{\text{$S$-gr.}}(\sG))$. By \cite{SGA3IIInew}*{XXIV, 1.3}, there is an exact sequence of \'etale $S$-sheaves:
\[
1\ra \sG/\un{\mathrm{Centr}}(\sG)\ra \underline{\mathrm{Aut}}(\sG)\ra \underline{\mathrm{Out}}(\sG)\ra 1,
\]
where $\un{\mathrm{Out}}(\sG)$ is the group scheme of outer automorphisms of $\sG$
and is represented by a locally constant group scheme, whose stalk at every geometric point of $S$ is a finitely generated group.
 We consider the following commutative diagram of pointed sets with exact rows
\[
\tikz {
\node  (A) at (-0.6,0) {$H^0(S, \un{\mathrm{Out}}(\sG))$};
\node  (B) at (3,0) {$H_{\et}^1(S, \sG/\mathrm{Centr}(\sG))$};
\node  (C) at (7,0) {$H_{\et}^1(S, \un{\mathrm{Aut}}(\sG))$};
\node  (D) at (11,0) {$H_{\et}^1(S, \un{\mathrm{Out}}(\sG))$};
\node  (a) at (-0.6,-1.5) {$H^0(K, \un{\mathrm{Out}}(\sG))$};
\node  (b) at (3,-1.5)   {$H_{\et}^1(K, \sG/\mathrm{Centr}(\sG))$};
\node  (c) at (7,-1.5) {$H_{\et}^1(K, \un{\mathrm{Aut}}(\sG))$};
\node  (d) at (11,-1.5) {$H_{\et}^1(K, \un{\mathrm{Out}}(\sG))$,};
\node  (l) at (3.3,-0.75)  {$\scriptstyle{f_1}$};
\node  (ll) at (7.3,-0.75) {$\scriptstyle{f_2}$};
\node  (lll) at (11.3,-0.75) {$\scriptstyle{f_3}$};
\node  (llll) at (-0.3,-0.75)  {$\scriptstyle{f_0}$};
\draw [draw = black, thin,
arrows={
- Stealth }]
(A) edge  (B);
\draw [draw = black, thin,
arrows={
- Stealth }]
(A) edge  (a);
\draw [draw = black, thin,
arrows={
- Stealth }]
(B) edge  (C);
\draw [draw = black, thin,
arrows={
- Stealth }]
(C) edge (D);
\draw [draw = black, thin,
arrows={
- Stealth }]
(a) edge (b);
\draw [draw = black, thin,
arrows={
- Stealth }]
(b) edge (c);
\draw [draw = black, thin,
arrows={
- Stealth }]
(c) edge (d);
\draw [draw = black, thin,
arrows={
- Stealth }]
(B) edge (b);
\draw [draw = black, thin,
arrows={
- Stealth }]
(C) edge (c);
\draw [draw = black, thin,
arrows={
- Stealth }]
(D) edge (d);
}
\]
where $f_2(\alpha)=0$. By assumption, the map $f_1$ is injective. In order to show that $\alpha=0$, we prove that the kernel of $f_3$ is trivial. Let $\cT$ be a $\outauto(\sG)$-torsor over $S$ that becomes trivial over $K$.
Then $\cT(K)\neq \emptyset$ contains a section $t_0$. For a surjective \'etale covering $X\pr\ra \Spec S$ such that $\cT|X\pr$ is a constant sheaf, the base change $X\pr_K$ satisfies $\cT(X\pr_K)=\cT(X\pr)$ and $\cT(X_K\pr\times_K X_K\pr)=\cT(X\pr\times_{\Spec S} X\pr)$, which fit into the following commutative diagram with exact rows:
\[
\tikz {
\node  (A) at (0,0) {$\cT(K)$};
\node  (B) at (2,0) {$\cT(X\pr_K)$};
\node  (C) at (5,0) {$\cT(X\pr_K\times_K X\pr_K)$};
\node  (a) at (0,-1.3)   {$\cT(S)$};
\node  (b) at (2,-1.3) {$\cT(X\pr)$};
\node  (c) at (5.2,-1.3) {$\cT(X\pr\times_{\Spec S} X\pr)$,};
\draw [draw = black, thin,
arrows={
- Stealth }]
(A) edge  (B);
\draw [draw = black, thin,
arrows={
- Stealth }]
(a) edge (b);
\draw [draw = black, thin,
arrows={
Hooks[right]- Stealth }]
(a) edge (A);
\draw [ - Stealth] (2.7cm,-0.2em) -- (3.7cm,-0.2em);
\draw [ - Stealth] (2.7cm,0.2em) -- (3.7cm,0.2em);
\draw [ - Stealth] (2.7cm,-1.23cm) -- (3.7cm,-1.23cm);
\draw [ - Stealth] (2.7cm,-1.37cm) -- (3.7cm,-1.37cm);
\draw (5.04cm,-0.3cm) -- (5.04cm,-1.0cm);
\draw (4.96cm,-0.3cm) -- (4.96cm,-1.0cm);
\draw (2.04cm,-0.3cm) -- (2.04cm,-1.0cm);
\draw (1.96cm,-0.3cm) -- (1.96cm,-1.0cm);
\draw [draw = black, thin,
arrows={
- Stealth }]
}
\]
such that the two images of $t_0$ under the double arrows in $\cT(X\pr)$ coincide in $\cT(X\pr\times_{\Spec S} X\pr)$. 
It follows that $\cT(S)=\cT(K)\neq \emptyset$ and $f_{3}$ has trivial kernel.
This argument also shows that $f_{0}$ is an isomorphism.
By diagram chasing, we conclude that $f_2$ has trivial kernel and $\sG\isom \sG'$.
\epf
In fact, the argument above also proves the following corollary.
\bcor
For a reductive group scheme $G$ over a regular semilocal ring $S$ with total ring of fractions $K$, if every form $G\pr$ of $G/\mathrm{Centr}(G)$ satisfies the injectivity of $H^1_{\et}(S, G\pr)\ra H^1_{\et}(K, G\pr)$, then
\[
\text{$H^1_{\et}(S, \un{\mathrm{Aut}}(G))\ra H^1_{\et}(K,\un{\mathrm{Aut}}(G))$\qq is injective.}
\]
\ecor
By \Cref{GSconj} and the main theorems in \cite{FP15, Pan17}, we obtain the following corollaries.
\bcor\label{modelsemilocalDedekind}
For a semilocal Dedekind ring $S$ with total ring of fractions $K$, any reductive $K$-group has at most one reductive $S$-model.
\ecor
\bcor
For a regular semilocal ring $S$ containing a field with total ring of fractions $K$, any reductive $K$-group has at most one reductive $S$-model.
\ecor


\bib{BF17}{article}{
  author={Bayer-Fluckiger, E.},
  author={First, U. A.},
  title={Rationally isomorphic hermitian forms and torsors of some non-reductive groups},
  journal={Adv. Math.},
  fjournal={Advances in Mathematics},
  volume={312},
  year={2017},
  pages={150--184},
  issn={0001-8708},
  doi={10.1016/j.aim.2017.03.012},
  url={https://doi-org.revues.math.u-psud.fr/10.1016/j.aim.2017.03.012},
}

\bib{BF18}{article}{
  author={Bayer-Fluckiger, E.},
  author={First, U. A.},
  author={Huruguen, M.},
  title={Orders that are \'{e}tale-locally isomorphic},
  journal={Algebra i Analiz},
  volume={31},
  year={2019},
  number={4},
  pages={1--15},
  issn={0234-0852},
}

\bib{BVG14}{article}{
  author={Beke, Sofie},
  author={Van Geel, Jan},
  title={An isomorphism problem for {A}zumaya algebras with involution over semilocal {B}\'{e}zout domains},
  journal={Algebr. Represent. Theory},
  fjournal={Algebras and Representation Theory},
  volume={17},
  year={2014},
  number={6},
  pages={1635--1655},
  issn={1386-923X},
  doi={10.1007/s10468-013-9463-6},
  url={https://doi-org.revues.math.u-psud.fr/10.1007/s10468-013-9463-6},
}

\bib{BB70}{article}{
  author={Bia\l ynicki-Birula, A.},
  title={Rationally trivial homogeneous principal fibrations of schemes},
  journal={Invent. Math.},
  volume={11},
  year={1970},
  pages={259--262},
  issn={0020-9910},
  doi={10.1007/BF01404652},
  url={https://doi-org.revues.math.u-psud.fr/10.1007/BF01404652},
}

\bib{BrT1}{article}{
  author={Bruhat, F.},
  author={J. Tits},
  title={Groupes r\'eductifs sur un corps local. {I}. {D}onn\'ees radicielles valu\'ees},
  journal={Inst. Hautes \'Etudes Sci. Publ. Math.},
  year={1972},
  number={41},
  pages={5-251},
  issn={0073-8301},
  url={http://www.numdam.org.revues.math.u-psud.fr:2048/item?id=PMIHES_1972__41__5_0},
  label={BT$_{\text {I}}$},
}

\bib{BrT2}{article}{
  author={Bruhat, F.},
  author={J. Tits},
  title={Groupes r\'eductifs sur un corps local. {II}. {S}ch\'emas en groupes. {E}xistence d'une donn\'ee radicielle valu\'ee},
  journal={Inst. Hautes \'Etudes Sci. Publ. Math.},
  year={1984},
  number={60},
  pages={197-376},
  issn={0073-8301},
  url={http://www.numdam.org.revues.math.u-psud.fr:2048/item?id=PMIHES_1984__60__5_0},
  label={BT$_{\text {II}}$},
}

\bib{BrT3}{article}{
  author={Bruhat, F.},
  author={J. Tits},
  title={Groupes alg\'{e}briques sur un corps local. {C}hapitre {III}. {C}ompl\'{e}ments et applications \`a la cohomologie galoisienne},
  journal={J. Fac. Sci. Univ. Tokyo Sect. IA Math.},
  volume={34},
  year={1987},
  number={3},
  pages={671--698},
  issn={0040-8980},
  label={BT$_{\text {III}}$},
}

\bib{Ces15d}{article}{
  author={{\v {C}}esnavi{\v {c}}ius, K{\k {e}}stutis},
  title={Topology on cohomology of local fields},
  journal={Forum Math. Sigma},
  volume={3},
  date={2015},
  pages={e16, 55},
  issn={2050-5094},
  review={\MR {3482265}},
  doi={10.1017/fms.2015.18},
}

\bib{CTS79}{article}{
  author={Colliot-Th\'el\`ene, Jean-Louis},
  author={Sansuc, Jean-Jacques},
  title={Fibr\'{e}s quadratiques et composantes connexes r\'{e}elles},
  journal={Math. Ann.},
  fjournal={Mathematische Annalen},
  volume={244},
  year={1979},
  number={2},
  pages={105--134},
  issn={0025-5831},
  doi={10.1007/BF01420486},
  url={https://doi.org/10.1007/BF01420486},
}

\bib{CTS87}{article}{
  author={Colliot-Th\'el\`ene, Jean-Louis},
  author={Sansuc, Jean-Jacques},
  title={Principal homogeneous spaces under flasque tori: applications},
  journal={J. Algebra},
  volume={106},
  date={1987},
  number={1},
  pages={148--205},
  issn={0021-8693},
  review={\MR {878473}},
  doi={10.1016/0021-8693(87)90026-3},
}

\bib{CTO92}{article}{
  author={Colliot-Th\'{e}l\`ene, Jean-Louis},
  author={Ojanguren, Manuel},
  title={Espaces principaux homog\`enes localement triviaux},
  journal={Inst. Hautes \'{E}tudes Sci. Publ. Math.},
  number={75},
  year={1992},
  pages={97--122},
  issn={0073-8301},
  url={http://www.numdam.org.revues.math.u-psud.fr:2048/item?id=PMIHES_1992__75__97_0},
}

\bib{Con12}{article}{
  author={Conrad, Brian},
  title={Weil and {G}rothendieck approaches to adelic points},
  journal={Enseign. Math. (2)},
  fjournal={L'Enseignement Math\'{e}matique. Revue Internationale. 2e S\'{e}rie},
  volume={58},
  year={2012},
  number={1-2},
  pages={61--97},
  issn={0013-8584},
  doi={10.4171/LEM/58-1-3},
  url={https://doi-org.revues.math.u-psud.fr/10.4171/LEM/58-1-3},
}

\bib{FR18}{article}{
  title={On the Grothendieck--Serre Conjecture about principal bundles and its generalizations},
  author={Fedorov, Roman},
  journal={arXiv preprint arXiv:1810.11844},
  year={2018},
}

\bib{FP15}{article}{
  author={Fedorov, Roman},
  author={Panin, Ivan},
  title={A proof of the Grothendieck--Serre conjecture on principal bundles over regular local rings containing infinite fields},
  journal={Publ. Math. Inst. Hautes \'Etudes Sci.},
  volume={122},
  date={2015},
  pages={169--193},
  issn={0073-8301},
  review={\MR {3415067}},
  doi={10.1007/s10240-015-0075-z},
}

\bib{Gir71}{article}{
  author={Giraud, Jean},
  title={Cohomologie non ab\'{e}lienne},
  note={Die Grundlehren der mathematischen Wissenschaften, Band 179},
  publisher={Springer-Verlag, Berlin-New York},
  year={1971},
  pages={ix+467},
}

\bib{Gro58}{article}{
  title={Torsion homologique et sections rationnelles},
  author={Grothendieck, Alexander},
  journal={Anneaux de Chow et applications, S{\'e}minaire Chevalley, 2e ann{\'e}e, Secr{\'e}tariat math{\'e}matique, Paris},
  year={1958},
}

\bib{Gro68}{article}{
  author={Grothendieck, Alexander},
  title={Le groupe de {B}rauer. {II}. {T}h\'{e}orie cohomologique},
  booktitle={Dix expos\'{e}s sur la cohomologie des sch\'{e}mas},
  series={Adv. Stud. Pure Math.},
  volume={3},
  pages={67--87},
  publisher={North-Holland, Amsterdam},
  year={1968},
}

\bib{Harder}{article}{
  author={Harder, G\"{u}nter},
  title={Eine {B}emerkung zum schwachen {A}pproximationssatz},
  journal={Arch. Math. (Basel)},
  volume={19},
  year={1968},
  pages={465--471},
  issn={0003-889X},
}

\bib{MB14}{article}{
  author={Moret-Bailly, Laurent},
  title={Un probl\`eme de descente},
  journal={Bull. Soc. Math. France},
  volume={124},
  year={1996},
  number={4},
  pages={559--585},
  issn={0037-9484},
  url={http://www.numdam.org.revues.math.u-psud.fr:2048/item?id=BSMF_1996__124_4_559_0},
}

\bib{nisnevich1983etale}{book}{
  title={ETALE COHOMOLOGY AND ARITHMETIC OF SEMISIMPLE GROUPS.},
  author={Nisnevich, Yevsey A.},
  year={1982},
  note={Thesis (Ph.D.)--Harvard University},
  publisher={ProQuest LLC, Ann Arbor, MI},
  pages={107},
  url={http://gateway.proquest.com/openurl?url_ver=Z39.88-2004&rft_val_fmt=info:ofi/fmt:kev:mtx:dissertation&res_dat=xri:pqdiss&rft_dat=xri:pqdiss:8303451},
}

\bib{Nis84}{article}{
  author={Nisnevich, Yevsey A.},
  title={Espaces homog\`enes principaux rationnellement triviaux et arithm\'etique des sch\'emas en groupes r\'eductifs sur les anneaux de Dedekind},
  language={French, with English summary},
  journal={C. R. Acad. Sci. Paris S\'er. I Math.},
  volume={299},
  date={1984},
  number={1},
  pages={5--8},
  issn={0249-6291},
  review={\MR {756297}},
}

\bib{Nis89}{article}{
  author={Nisnevich, Yevsey A.},
  title={Rationally trivial principal homogeneous spaces, purity and arithmetic of reductive group schemes over extensions of two-dimensional regular local rings},
  language={English, with French summary},
  journal={C. R. Acad. Sci. Paris S\'er. I Math.},
  volume={309},
  date={1989},
  number={10},
  pages={651--655},
  issn={0764-4442},
  review={\MR {1054270}},
}

\bib{Pan15}{article}{
  author={Panin, Ivan Alexandrovich},
  title={Proof of Grothendieck--Serre conjecture on principal bundles over regular local rings containing a finite field},
  journal={Linear algebraic groups and Related Structures},
  volume={559},
  date={2015},
  pages={42 pages},
  ulr={https://www.math.uni-bielefeld.de/lag/},
}

\bib{Pan17}{unpublished}{
  author={Panin, Ivan},
  title={Proof of {G}rothendieck--Serre's conjecture on principal bundles over regular local rings containing a finite field},
  note={preprint, avaiable at https://arxiv.org/abs/1707.01767},
  year={2017},
}

\bib{PS17}{article}{
  author={Panin, Ivan},
  author={Stavrova, A. K.},
  title={On the {G}rothendieck--{S}erre conjecture concerning principal {G}-bundles over semi-local {D}edekind domains},
  note={Reprinted in J. Math. Sci. (N.Y.) {{\bf {2}}22}, (2017), no. 4, 453--462},
  journal={Zap. Nauchn. Sem. S.-Peterburg. Otdel. Mat. Inst. Steklov. (POMI)},
  volume={443},
  year={2016},
  number={Voprosy Teorii Predstavleni\u {\i } Algebr i Grupp. 29},
  pages={133--146},
  issn={0373-2703},
  doi={10.1007/s10958-017-3316-5},
  url={https://doi-org.revues.math.u-psud.fr/10.1007/s10958-017-3316-5},
}

\bib{PV}{book}{
  author={Platonov, Vladimir},
  author={Rapinchuk, Andrei},
  title={Algebraic groups and number theory},
  series={Pure and Applied Mathematics},
  volume={139},
  note={Translated from the 1991 Russian original by Rachel Rowen},
  publisher={Academic Press, Inc., Boston, MA},
  year={1994},
  pages={xii+614},
  isbn={0-12-558180-7},
}

\bib{PY06}{article}{
  author={Prasad, Gopal},
  author={Yu, Jiu-Kang},
  title={On quasi-reductive group schemes},
  note={With an appendix by Brian Conrad},
  journal={J. Algebraic Geom.},
  volume={15},
  year={2006},
  number={3},
  pages={507--549},
  issn={1056-3911},
  doi={10.1090/S1056-3911-06-00422-X},
  url={https://doi-org.revues.math.u-psud.fr/10.1090/S1056-3911-06-00422-X},
}

\bib{Sch13}{article}{
  author={Scholze, Peter},
  title={{$p$}-adic {H}odge theory for rigid-analytic varieties},
  journal={Forum Math. Pi},
  volume={1},
  year={2013},
  pages={e1, 77},
  issn={2050-5086},
  doi={10.1017/fmp.2013.1},
  url={https://doi-org.revues.math.u-psud.fr/10.1017/fmp.2013.1},
}

\bib{Ser58}{article}{
  title={Espaces fibr{\'e}s alg{\'e}briques},
  author={Serre, Jean-Pierre},
  journal={S{\'e}minaire Claude Chevalley},
  volume={3},
  pages={1--37},
  year={1958},
}

\bib{Ser02}{book}{
  author={Serre, Jean-Pierre},
  title={Galois cohomology},
  series={Springer Monographs in Mathematics},
  edition={Corrected reprint of the 1997 English edition},
  note={Translated from the French by Patrick Ion and revised by the author},
  publisher={Springer-Verlag},
  place={Berlin},
  date={2002},
  pages={x+210},
  isbn={3-540-42192-0},
  review={\MR {1867431 (2002i:12004)}},
}

\bib{SGA3II}{book}{
  title={Sch\'emas en groupes. II: Groupes de type multiplicatif, et structure des sch\'emas en groupes g\'en\'eraux},
  language={French},
  series={S\'eminaire de G\'eom\'etrie Alg\'ebrique du Bois Marie 1962/64 (SGA 3). Dirig\'e par M. Demazure et A. Grothendieck. Lecture Notes in Mathematics, Vol. 152},
  publisher={Springer-Verlag, Berlin-New York},
  date={1970},
  pages={ix+654},
  review={\MR {0274459 (43 \#223b)}},
  label={SGA~3$_{\text {II}}$},
}

\bib{SGA3IIInew}{collection}{
  title={Sch\'emas en groupes (SGA 3). Tome III. Structure des sch\'emas en groupes r\'eductifs},
  language={French},
  series={Documents Math\'ematiques (Paris) [Mathematical Documents (Paris)], 8},
  editor={Gille, Philippe},
  editor={Polo, Patrick},
  note={S\'eminaire de G\'eom\'etrie Alg\'ebrique du Bois Marie 1962--64. [Algebraic Geometry Seminar of Bois Marie 1962--64]; A seminar directed by M. Demazure and A. Grothendieck with the collaboration of M. Artin, J.-E. Bertin, P. Gabriel, M. Raynaud and J-P. Serre; Revised and annotated edition of the 1970 French original},
  publisher={Soci\'et\'e Math\'ematique de France, Paris},
  date={2011},
  pages={lvi+337},
  isbn={978-2-85629-324-9},
  review={\MR {2867622}},
  label={SGA~3$_{\text {III new}}$},
}

\bib{SGA4II}{book}{
  title={Th\'eorie des topos et cohomologie \'etale des sch\'emas. Tome 2},
  language={French},
  series={Lecture Notes in Mathematics, Vol. 270},
  note={S\'eminaire de G\'eom\'etrie Alg\'ebrique du Bois-Marie 1963--1964 (SGA 4); Dirig\'e par M. Artin, A. Grothendieck et J. L. Verdier. Avec la collaboration de N. Bourbaki, P. Deligne et B. Saint-Donat},
  publisher={Springer-Verlag},
  place={Berlin},
  date={1972},
  pages={iv+418},
  review={\MR {0354653 (50 \#7131)}},
  label={SGA~4$_{\text {II}}$},
}

\bib{SP}{misc}{
  shorthand={Stacks},
  author={The {Stacks Project Authors}},
  title={\textit {Stacks Project}},
  howpublished={\url {https://stacks.math.columbia.edu}},
  year={2020},
  label={SP},
}

\bib{JKY}{incollection}{
  author={Yu, Jiu-Kang},
  title={Bruhat--{T}its theory and buildings},
  booktitle={Ottawa lectures on admissible representations of reductive {$p$}-adic groups},
  series={Fields Inst. Monogr.},
  volume={26},
  pages={53--77},
  publisher={Amer. Math. Soc., Providence, RI},
  year={2009},
}



\begin{bibdiv}
\begin{biblist}

\bibselect{bibliography}
\bibliographystyle{alpha}
\end{biblist}
\end{bibdiv}

\end{document}